\documentclass[article]{aiaa-pretty}    

\usepackage[utf8]{inputenc}  

\usepackage{amssymb}
\usepackage{amsmath}
\usepackage{booktabs}
 \usepackage{varioref}
 \usepackage{wrapfig}
 \usepackage{threeparttable}
 \usepackage{dcolumn}
  \newcolumntype{d}{D{.}{.}{-1}}
 \usepackage{nomencl}
  \makeglossary
 \usepackage{subfigure}
 \usepackage{graphicx}
 \usepackage{subfigmat}
 \usepackage{fancyvrb}
  \fvset{fontsize=\footnotesize,xleftmargin=2em}
 \usepackage{lettrine}
\usepackage{cite}

\newcommand{\Real}{\mathbb R}
\newcommand{\set}[1]{\left\{#1\right\}}
\newcommand{\real}[1]{{\mathbb R}^{#1}}
\newcommand{\be}{{\boldsymbol e}}
\newcommand{\bff}{{\boldsymbol f}}
\newcommand{\bg}{{\boldsymbol g}}
\newcommand{\bh}{{\boldsymbol h}}

\newcommand{\bp}{{\boldsymbol p}}
\newcommand{\bq}{{\boldsymbol q}}

\newcommand{\bu}{{\boldsymbol u}}

\newcommand{\bx}{\boldsymbol x}

\newcommand{\bxf}{{\bx(\cdot)}}  

\newcommand{\buf}{{\bu(\cdot)}}  

\newcommand{\bP}{{\boldsymbol P}}

\newcommand{\U}{\mathbb{U}}

\newcommand{\bzero}{{\bf 0}}

\newcommand{\bX}{{\mbox{\boldmath $X$}}}

\newcommand{\blam}{{\mbox{\boldmath $\lambda$}}}

\newcommand{\bomega}{\mbox{\boldmath$\omega$}}

\newcommand{\bsigma}{\mbox{\boldmath$\sigma$}}
\newcommand{\bW}{{\mbox{\boldmath $W$}}}
\newcommand{\bLam}{\mbox{\boldmath$\Lambda$}}

\author{
I. M. Ross,\thanks{Distinguished Professor and Program Director, Control and Optimization, Department of Mechanical and Aerospace Engineering}
R. J. Proulx\thanks{Research Professor, Control and Optimization Laboratories, Space Systems Academic Group}
and
M.~Karpenko\thanks{Research Professor and Associate Director, Control and Optimization Laboratories, Department of Mechanical and Aerospace Engineering}
\\
\textit{Naval Postgraduate School, Monterey, CA 93943}
}

\title{Unscented Trajectory Optimization}

\abstract{
In a nutshell, unscented trajectory optimization is the generation of optimal trajectories through the use of an unscented transform.  Although unscented trajectory optimization was introduced by the authors about a decade ago, it is reintroduced in this paper as a special instantiation of tychastic optimal control theory. Tychastic optimal control theory (from \textit{Tyche}, the Greek goddess of chance) avoids the use of a Brownian motion and the resulting It\^{o} calculus even though it uses random variables across the entire spectrum of a problem formulation.  This approach circumvents the enormous technical and numerical challenges associated with stochastic trajectory optimization.  Furthermore, it is shown how a tychastic optimal control problem that involves nonlinear transformations of the expectation operator can be quickly instantiated using an unscented transform. These nonlinear transformations are particularly useful in managing trajectory dispersions be it associated with path constraints or targeted values of final-time conditions.
This paper also presents a systematic and rapid process for formulating and computing the most desirable tychastic trajectory using an unscented transform. Numerical examples are used to illustrate how unscented trajectory optimization may be used for risk reduction and mission recovery caused by uncertainties and failures.

}		

\begin{document}

\maketitle{}

\section{Introduction}

\begin{figure}[!ht]
	\centering
    \includegraphics[angle = 0, width=0.65\columnwidth,clip]{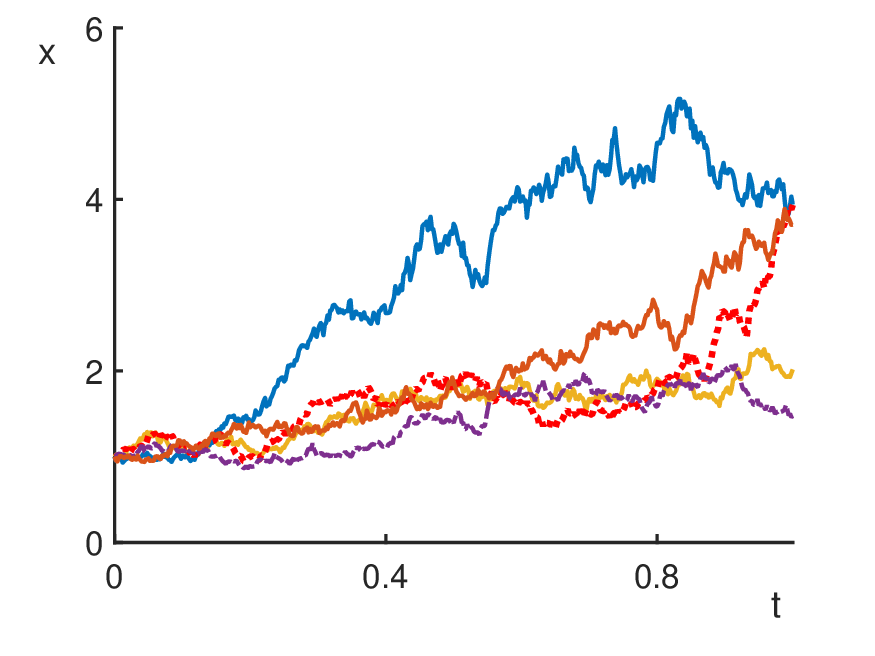}
	\caption{An illustration of the evolutions of state trajectories for a stochastic differential equation given by \eqref{eq:SDE-Ito}.}
	\label{fig:stocTrajs}
\end{figure}

Practical trajectory optimization problems naturally contain various uncertainties.  The conventional mathematical model for a system with uncertainties is a \emph{\textbf{stochastic control differential equation}} given by\cite{kushner,peng,RMP-book,yong-zhou-1999Book},
\begin{equation}\label{eq:SCDE-Ito}
d\bx = \bff(\bx, \bu, t)dt + \bsigma(\bx, \bu, t)\,d\bW
\end{equation}
where, $\bx \in \real{N_x}$ is an $N_x$-dimensional state variable, $\bW(\cdot)$ is a standard $N_w$-dimensional Wiener process (i.e., Brownian motion), $\bsigma$ is an $N_x \times N_w$ diffusion matrix that may depend on the control variable $\bu \in \real{N_u}$ (in addition to $\bx$ and $t$) and $\bff$ is a deterministic dynamics function.  Suppose a deterministic control trajectory, $ t \mapsto \bu^\sharp(t)$, and an initial condition, $\bx(0) = \bx^0$, are given; then, \eqref{eq:SCDE-Ito} becomes a stochastic ordinary differential equation of the It\^{o}-type \cite{kushner,peng},
\begin{equation}\label{eq:SDE-Ito}
d\bx = \bff(\bx, \bu^\sharp(t), t)dt + \bsigma(\bx, \bu^\sharp(t), t)\,d\bW
\end{equation}
Let $\bx^\sharp(t)$ be a solution to \eqref{eq:SDE-Ito}. As illustrated in Fig.~\ref{fig:stocTrajs}, $\bx^\sharp(t)$ is given by\cite{SDE-sim-book,SDE-num-2001},
\begin{equation}\label{eq:SDE-sol}
\bx^\sharp(t) = \bx^0 + \int_0^t \bff(\bx(s), \bu^\sharp(s), s)\, ds + \int_0^t \bsigma(\bx(s), \bu^\sharp(s), s)\,d\bW
\end{equation}
where, the stochastic integral, $\int_0^t (\cdot)\,d\bW $ must be computed using the rules of It\^{o} calculus\cite{kushner,peng,yong-zhou-1999Book}.
%
%
For example, the stochastic Euler method for \eqref{eq:SDE-Ito}, known as the Euler-Maruyama method\cite{SDE-sim-book,SDE-num-2001}, is given by,
\begin{equation}\label{eq:Euler-Maruyama}
\bx(t + h) = \bx(t) +  \bff(\bx(t), \bu^\sharp(t), t) h + \bsigma(\bx(t), \bu^\sharp(t), t)\bW_h \sqrt{h}
\end{equation}
where, $h$ is the step-size and $\bW_h \sim \mathcal{N}(\mathbf{0},\mathbf{1})$. Note the $\sqrt{h}$ in \eqref{eq:Euler-Maruyama}.
That is, simulating \eqref{eq:SDE-Ito} requires stochastic numerical methods\cite{SDE-num-2001} that are sharply different from deterministic numerical methods\cite{hnw-ode}. \emph{\textbf{Erroneous results are obtained\cite{SDE-sim-book} if \eqref{eq:SDE-Ito} is solved using standard Runge-Kutta methods}}.

What makes \eqref{eq:SCDE-Ito} even more challenging as a model for trajectory optimization is that the stochastic maximum principle\cite{peng,RMP-book, poznyak-2002} involves two Hamiltonians and two adjoint equations, neither of which resemble the deterministic Pontryagin equations\cite{ross-book}. For example, the first Hamiltonian for \eqref{eq:SCDE-Ito} is given by the sum of a Pontryagin-type Hamiltonian and the trace of a term involving the diffusion matrix\cite{peng,poznyak-2002,yong-zhou-1999Book},
\begin{equation}\label{eq:Hamil4stoc}
H(\blam, \bLam, \bx, \bu, t) := \blam^T\bff(\bx, \bu, t) + \text{tr}[\bLam^T\bsigma(\bx, \bu, t)]
\end{equation}
where, $\bLam \in \real{N_x \times N_w}$ is an adjoint diffusion matrix associated with an adjoint stochastic differential equation involving $\blam$.  It turns out that an optimal control does not even maximize $H(\blam, \bLam, \bx, \bu, t)$\cite{peng}; instead, it maximizes a second Hamiltonian that depends upon \eqref{eq:Hamil4stoc} and an additional new adjoint matrix that satisfies a second-order adjoint process\cite{peng,RMP-book,poznyak-2002}.  Furthermore, the second-order matrix adjoint equation involves the Hessian of the first Hamiltonian\cite{peng}. For the purposes of brevity we do not provide these lengthy details here.  Instead, we note the following: Although the mathematics of a stochastic optimal control problem is sufficiently advanced, it is not clear how to use it effectively in formulating and solving a stochastic trajectory optimization problem. For example, the diffusion term (see \eqref{eq:SDE-Ito} and \eqref{eq:Euler-Maruyama}) can induce a very high degree of sensitivity\cite{ross-book} if a shooting method\cite{shooting-2013,ross-book,trelat:survey} is used to solve the resulting stochastic trajectory optimization problem. Furthermore, the numerical iterations would not be repeatable due to the presence of the Wiener process. Similarly, because a standard Runge-Kutta method is inapplicable\cite{SDE-sim-book,SDE-num-2001} as a discretization method for  \eqref{eq:SCDE-Ito}, the same would hold for a collocation method\cite{conway:survey,trelat:survey}. \emph{\textbf{In fact, to the best of the authors' knowledge, the stochastic maximum principle\cite{kushner,peng,RMP-book,yong-zhou-1999Book} has never been used to address the formulation or optimality of a nontrivial, stochastic, nonlinear, aerospace trajectory optimization problem.}}

\begin{figure}[!ht]
	\centering
    \includegraphics[angle = 0, width=0.65\columnwidth,clip]{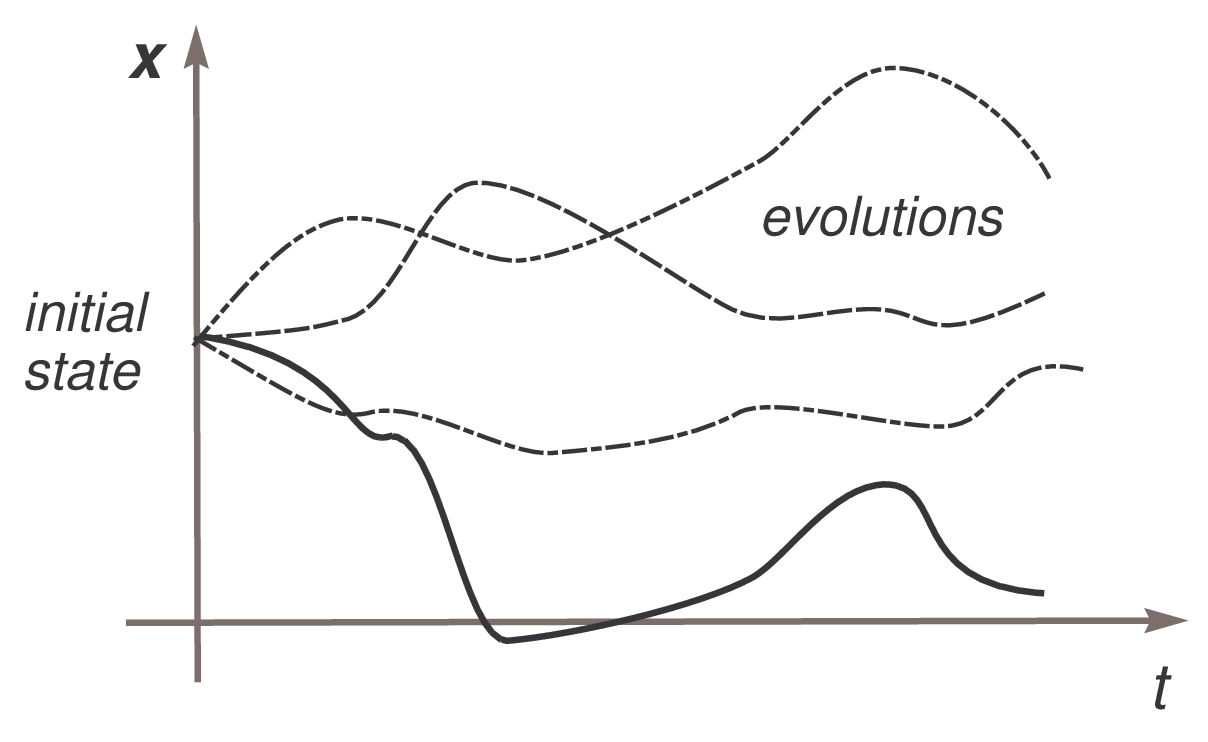}
	\caption{Schematic for the evolutions of state trajectories for a tychastic differential equation; compare with Fig.~\ref{fig:stocTrajs}.  Figure adapted from [\citen{LS:acc-2}] and [\citen{aubin=tyche}].}
	\label{fig:DiffInc}
\end{figure}

From the arguments of the preceding paragraphs, it is clear that the main source of all of the mathematical and computational difficulties in a stochastic trajectory optimization problem is the use of the Wiener process for modeling randomness.  In circumventing the difficulties stemming from the use of the Wiener process,
we introduced in [\citen{uoc-1:issfd,uoc-2:SD,ug:acc,RS:jgcd,LS:acc-1,LS:acc-2,uo-1:2015,uo-1:MRYpts-2016, uo-patent-2016,uo-patent-2017,uoc-patent-1,uoc-patent-2}] a series of concepts and methods that culminated with the theory of \emph{\textbf{tychastic trajectory optimization}}\cite{LS:acc-2,ross-book} (from \textit{Tyche}, the Greek goddess of chance). In a tychastic trajectory optimization problem, the stochastic differential equation given by \eqref{eq:SCDE-Ito} is abandoned in favor of a tychastic differential equation\cite{ross-book},
\begin{equation}\label{eq:tyc-code}
\dot\bx = \bff(\bx, \bu, t; \bp), \quad \bp \in supp(\bp) \subseteq \real{N_p}
\end{equation}
where, $\bp$ is an $N_p$-dimensional random variable defined over a support $supp(\bp)$ that is not necessarily compact.  In their original definition, Aubin et al\cite{aubin=tyche} define tyches more generally than that implied in \eqref{eq:tyc-code}.  In this paper, we use our simpler form while noting that the presence of $\bp$ is not necessarily exclusive to uncertainties in the dynamics but also in the initial conditions, path constraints and the cost functional itself\cite{LS:acc-1,LS:acc-2,uoc-patent-1,uoc-patent-2,ross-book}.  Deferring a discussion of these generalities, we first note that the dynamics function in \eqref{eq:tyc-code} is conditionally deterministic; i.e., $\bff$ is deterministic if $\bp$ is known.  Furthermore, unlike \eqref{eq:SCDE-Ito}, the derivative in \eqref{eq:tyc-code} is the ordinary derivative. Because $\bp$ is indeed unknown and uncertain, we view \eqref{eq:tyc-code} in terms of a deterministic selection of a controlled differential inclusion\cite{clsw},
\begin{equation}\label{eq:cdi}
\dot\bx \in \mathcal{F}(\bx, \bu, t) := \set{\bff(\bx, \bu, t; \bp):\ \bp \in supp(\bp)}
\end{equation}
Thus, the evolution of \eqref{eq:tyc-code} for a deterministic control function $t \mapsto \bu^\sharp(t)$ is considered to be set-valued (see Fig.~\ref{fig:DiffInc}) and does not require the use of any exotic It\^{o} calculus.
%
Compare Fig.~\ref{fig:DiffInc} with Fig.~\ref{fig:stocTrajs}.
Such dynamical systems have found widespread applications in search theory\cite{koopman-2,stone,phelps-2:cdc} and quantum control\cite{Li:phys,Li:pnas}.  Equation~\eqref{eq:cdi} has also been widely studied in the Russian literature\cite{RussiaReport} and has led to a minimax or robust optimal control theory as put forth separately by Vinter\cite{vinter-minimax} and Boltyanski and Poznyak~\cite{RMP-book}.  In order to separate \eqref{eq:tyc-code} from \eqref{eq:SCDE-Ito} and the ensuing mathematics, it is called a tychastic differential equation\cite{aubin-alm-2012,ross-book}. Thus, in sharp contrast to a stochastic differential equation, we may now indeed use standard Runge-Kutta methods\cite{hnw-ode} to simulate \eqref{eq:tyc-code} for a given random value of $\bp \in supp(\bp)$ to produce trajectories such as the ones shown in Fig.~\ref{fig:DiffInc}.  Furthermore, we may select a sufficiently large number of random samples $N_s \in \mathbb{N}$ from $supp(\bp)$ to numerically construct the approximate statistics associated with $\mathcal{F}(\bx(t), \bu(t), t)$. Although it circumvents the myriad of challenges associated with stochastic optimal control theory, the computational challenge with such a Monte Carlo approach is that the generation of a single tychastic optimal control requires the equivalent of a solution to an ensemble of $N_s$ deterministic dynamical systems, where $N_s$ may be as large as $1,000$ or even $10,000$. Thus, for example, a Monte Carlo approach would require one to solve $1,000$--$10,000$ optimal control problems simultaneously. Unscented optimal control theory\cite{RS:jgcd,LS:acc-2,uoc-patent-1,uoc-patent-2} addresses this scale challenge head on: instead of selecting a large number of random samples from $supp(\bp)$, a  significantly smaller set of sigma points\cite{julier:simplex,julier-acc-95} is used for faster tychastic trajectory optimization.  That is, unscented trajectory optimization is indeed a fast instantiation of tychastic optimal control theory.

As implied earlier, different aspects of the preceding ideas have been presented elsewhere\cite{uoc-1:issfd,uoc-2:SD,uoc-patent-1,uoc-patent-2,ug:acc,RS:jgcd,LS:acc-1,LS:acc-2}.  In this paper, we present a new perspective that incorporates the concepts discussed in \cite{uoc-1:issfd,uoc-2:SD,ug:acc} via the overarching framework of tychastic trajectory optimization.  Tychastic trajectory optimization subsumes the Lebesgue-Stieltjes framework presented in \cite{LS:acc-1,LS:acc-2} which, in turn, generalized the Riemann-Stieltjes theory presented in \cite{RS:jgcd}. In the next section, we provide new and additional details on all these concepts using \eqref{eq:cdi} as a starting point.

\section{Tychastic Trajectory Optimization}
There are multiple ways to formulate a tychastic trajectory optimization problem\cite{ross-book}.  Here we discuss several non-exhaustive ways to formulate a problem so that the ``most desirable'' problem can be properly framed based on the need to control and optimize the relevant statistics associated with a given variable/function.

\subsection{Mathematical Preliminaries}
From \eqref{eq:cdi}, it follows that, for any given control trajectory $t \mapsto \bu(t)$, the evolution of a tychastic state trajectory may be written in terms of a collection of all possibilities (see Fig.~\ref{fig:DiffInc}) according to,
\begin{equation}
t \mapsto \set{\bx(t, \bp):\ \bp \in supp(\bp)}
\end{equation}
If one were to simply replace $\bxf$ in a conventional deterministic cost functional\cite{longuski,bryson:ho} $J: (\bxf, \buf, t_0, t_f) \mapsto \Real$ by $\bx(\cdot, \cdot)$, then $J$ itself may become uncertain\cite{ross-book} if the map,
\begin{equation}\label{eq:unc-map4J}
(\bx(\cdot, \cdot), \buf, t_0, t_f) \mapsto \Real
\end{equation}
depends explicitly on $\bx(\cdot, \cdot)$.  Hence, it may not be possible to minimize a conventional cost functional by simply replacing $\bxf$ by $\bx(\cdot, \cdot)$\cite{LS:acc-1,LS:acc-2,ross-book}.  However, because \eqref{eq:unc-map4J} may be construed as a deterministic cost functional for a given value of $\bp$, we can construct a new \textbf{\emph{functional of a functional}}\cite{ross-book} according to,
\begin{equation}\label{eq:tyc-cost}
J_{tyc}:\ \big(\bx(\cdot, \cdot), \buf, t_0, t_f; supp(\bp)\big) \mapsto \Real
\end{equation}
where, $J_{tyc}$ is a deterministic functional of a tychastic function $\bx(\cdot, \cdot)$, a tychastic/random parameter $\bp \in supp(\bp)$, a deterministic function $\buf$ and deterministic clock times $t_0$ and $t_f$.

\subsection{Sample Tychastic Cost Functionals}

Unlike a benchmark ``Bolza'' cost functional\cite{longuski,bryson:ho,clsw} that is widely used in deterministic optimal control theory, there is no standard computational formula for constructing \eqref{eq:tyc-cost}.  As a result, we define some useful tychastic cost functionals.

\subsubsection{Robust Cost Functional}
The authors of [\citen{RMP-book}] and [\citen{vinter-minimax}] propose a minimax or ``robust'' cost functional in terms of minimizing the maximal cost. In this case, the computational formula for \eqref{eq:tyc-cost} is the the worst-case tychastic cost functional given by,
\begin{equation}\label{eq:minimax}
J_{tyc}[\bx(\cdot, \cdot), \buf, t_0, t_f; supp(\bp)] :=  \mathop\text{ess\,sup}_{\bp \in supp(\bp)} E(\bx(t_f, \bp))
\end{equation}
where, $E$ is a ``randomized version'' of a given deterministic endpoint cost function, $E: \real{N_x} \to \Real$; i.e., a ``Mayer'' type\cite{longuski,ross-book,bryson:ho}, that is usually designed by considering deterministic situations.

\subsubsection{Average Cost Functional}

As has been noted elsewhere\cite{LS:acc-1,LS:acc-2,ross-book}, a consideration of \eqref{eq:minimax} as a standard bearer for tychastic optimal control is quite limiting.  In fact, there is a practical need to consider more general cost functionals including new ways of formulating uncertain endpoint conditions and path constraints all of which are not considered in [\citen{RMP-book}] and [\citen{vinter-minimax}].  In view of this, we begin by considering a tychastic cost functional that is different from \eqref{eq:minimax} and given by the ``average'' cost functional\cite{LS:acc-1,LS:acc-2,ross-book},
\begin{equation}\label{eq:tyc-avg}
J_{tyc}[(\bx(\cdot, \cdot), \buf, t_0, t_f; supp(\bp))] := 
\int \cdots \int_{supp(\bp)} J[(\bx(\cdot, \bp), \buf, t_0, t_f, \bp)]\, dm(\bp)
\end{equation}
where, $J$ is a given conditionally deterministic Bolza cost functional $(\bx(\cdot, \bp), \buf, t_0, t_f, \bp)$ $\mapsto \Real$ that is defined for a fixed value of $\bp \in supp(\bp)$, and $\int \cdots \int_{supp(\bp)}$ is an $N_p$-dimensional Lebesgue-Stieltjes integral with measure $m: \bp \mapsto \Real_+$.  Equation~\eqref{eq:tyc-avg} formed the basis of Riemann-Stieltjes optimal control theory\cite{RS:jgcd} that was subsequently generalized to a Lebesgue-Stieltjes formalism\cite{LS:acc-1,LS:acc-2}.

\subsubsection{A Dispersion Cost Functional}
It is clear that both \eqref{eq:minimax} and \eqref{eq:tyc-avg} are motivated by cost functionals from deterministic optimal control theory.
Consider next a tychastic cost functional given by,
\begin{equation}\label{eq:tyc-trCov}
J_{tyc}[\bx(\cdot, \cdot), \buf, t_0, t_f; supp(\bp)] := \text{tr}\,\text{Cov}[\bx(t_f, \bp)]
\end{equation}
where, tr Cov is the trace of the covariance (matrix) associated with $\bx(t_f, \bp)$ with Cov defined by,
\begin{equation}\label{eq:cov-def}
\text{Cov}[\bx(t_f, \bp)] := 
\mathcal{E}\left[ \big(\bx(t_f, \bp) - \mathcal{E}[\bx(t_f, \bp)] \big) \big(\bx(t_f, \bp) - \mathcal{E}[\bx(t_f, \bp)] \big)^T   \right]
\end{equation}
and $\mathcal{E}$ is the expectation operator given by,
\begin{equation}\label{eq:expect-def}
\mathcal{E}[\cdot]:= \int \cdots \int_{supp(\bp)} (\cdot)\, dm(\bp)
\end{equation}
It is apparent that neither of the computational formulas given by the right-hand-sides of \eqref{eq:minimax} or \eqref{eq:tyc-avg} can be used to compute \eqref{eq:tyc-trCov}.  Note also that \eqref{eq:minimax}, \eqref{eq:tyc-avg} and \eqref{eq:tyc-trCov} are all functionals of functionals -- a signature feature of tychastic optimal control theory\cite{ross-book}.  Instead of the trace of the covariance matrix, one can also choose its determinant or any other ``measure'' of the covariance matrix.

\subsubsection{Another Nonlinear Cost Functional}
Obviously, \eqref{eq:tyc-trCov} is a nonlinear cost functional (of a functional).  Yet another nonlinear cost functional may be written as,
\begin{equation}\label{eq:tyc-E-avg}
J_{tyc}[\bx(\cdot, \cdot), \buf, t_0, t_f; supp(\bp)] := 
E\big(\mathcal{E}[\bx(t_0, \bp)], \mathcal{E}[\bx(t_f, \bp)], t_0, t_f, \mathcal{E}[\bp]\big)
\end{equation}
where $E: \real{N_x} \times \real{N_x} \times \Real \times \Real \to \Real $ is a deterministic nonlinear endpoint cost.  Thus, $\mathcal{E}$ transforms nonlinearly according to $E$.  Furthermore, it is apparent that \eqref{eq:tyc-E-avg} cannot be computed using any of the formulas given by \eqref{eq:minimax}, \eqref{eq:tyc-avg} and \eqref{eq:tyc-trCov}.

\subsubsection{Event Probability Cost Functional}
As noted in \cite{ross-book}, one of the most difficult functionals to compute is given by,
\begin{equation}\label{eq:tyc-Pr}
J_{tyc}[\bx(\cdot, \cdot), \buf, t_0, t_f; supp(\bp)] :=  \text{Pr}\set{\be(\bx(t_f, \bp), t_f, \bp) \le \bzero}
\end{equation}
where $\text{Pr}\set{\cdot}$ denotes the probability of the event $\set{\cdot}$, and $\be$ is an event function\cite{ross-book}.  The difficulty in computing \eqref{eq:tyc-Pr} stems from it definition:
\begin{equation}\label{eq:Pr-def}
\text{Pr}\set{\be(\bx(t_f, \bp), t_f, \bp) \le \bzero} := \int\cdots\int_{\be(\bx(t_f, \bp), t_f, \bp) \le \bzero} \, dm(\bp)
\end{equation}
That is, unlike the expectation operator (see \eqref{eq:expect-def}) the domain of integration in \eqref{eq:Pr-def} is unknown a priori.

The preceding examples illustrate why a generic tychastic cost functional has not been found as in the case of deterministic optimal control theory\cite{ross-book,vinter,clsw,longuski,bryson:ho} wherein a Bolza functional is quite sufficient to model many types of performance indices.  Note also that the list of cost functionals defined in the preceding paragraphs are nonexhaustive.

\subsection{Risk, Reliability and Confidence Levels}
In a typical deterministic trajectory optimization problem\cite{ross-book,longuski,bryson:ho}, a target set $\mathcal{T}$ is specified in terms of function inequalities,
\begin{equation}\label{eq:T:=}
\mathcal{T} := \set{(\bx_f, t_f) \in \real{N_x} \times \Real:\ \be^L \le \be(\bx_f, t_f) \le \be^U }
\end{equation}
where $\be :\real{N_x} \times \Real \to \real{N_e} $ is a given endpoint/event function and $\be^L$ and $\be^U$ are the lower and upper bounds on the values of $\be$.  In a tychastic trajectory optimization problem it is quite possible (and very likely) that $\be(\bx(t_f, \bp), t_f)$ will violate one or more of its bounds. In this instance, one may impose a \emph{chance constraint}\cite{charnes-cooper,ross-book} given by,
\begin{equation}\label{eq:CC4T}
\text{Pr}\set{\bx(t_f, \bp) \in \mathcal{T}} := 
\int\cdots\int_{\bx(t_f, \bp) \in \mathcal{T}} \, dm(\bp)  \ge R_{\mathcal{T}} := 1- r_{\mathcal{T}}
\end{equation}
where,   $R_{\mathcal{T}} \in (0, 1)$ is a specified reliability or confidence level in achieving the target condition, $\bx(t_f, \bp) \in \mathcal{T}$. Alternatively, as implied in \eqref{eq:CC4T}, the chance constraint may specified in terms of the risk $r_{\mathcal{T}} \in (0,1)$ in not achieving the desired condition.

From \eqref{eq:CC4T} it follows that \eqref{eq:T:=} must be specified with the power of the continuum with non-zero, non-atomic measure in order to ensure that $\text{Pr}\set{\cdot} \ne 0$.  Surprisingly, the other extreme case of $\text{Pr}\set{\cdot} = 1$ is theoretically possible even when $\mathcal{T}$ is a singleton\cite{uoc-2:SD,ross-book}.  We denote this case as a \emph{transcendental target condition}\cite{ross-book,ug:acc} stipulated in the form of a semi-infinite (almost surely) constraint,
\begin{equation}\label{eq:efinal=trans}
\be^L \le \be(\bx(t_f, \bp), t_f, \bp) \le \be^U \quad \forall\ \bp \in supp(\bp)
\end{equation}
%

\subsection{Several Tychastic Problem Formulations}

In the same vein as \eqref{eq:efinal=trans}, a transcendental path constraint may be specified as,
\begin{equation}\label{eq:h=trans}
\bh^L \le \bh(\bx(t, \bp), \bu(t), t, \bp) \le \bh^U \quad \forall\ \bp \in supp(\bp)
\end{equation}
The central problem in using \eqref{eq:efinal=trans} and \eqref{eq:h=trans} to formulate generic boundary conditions and path constraints in a tychastic optimal control problem is that it is highly likely the resulting feasible set will be empty. Hence, a solution might not exist.  \emph{In fact, as noted before\cite{ug:acc,ross-book}, the existence of a solution is the main problem in tychastic optimal control.}  Nonetheless, an aspirational version of a family of transcendental tychastic optimal control problems may be framed in terms of a semi-infinite dimensional framework given by\cite{ug:acc,LS:acc-2}:
\begin{eqnarray}
& \forall\ \bp \in supp(\bp) & \nonumber\\[0.5em]
&(A^\infty) \left\{
\begin{array}{lll}
\text{Minimize }   &J_{tyc}[\bx(\cdot, \cdot), \buf, t_0, t_f; supp(\bp)] \\[0.5em]
\text{Subject to }  & \dot{\bx}(t, \bp) = {\bff}(\bx(t, \bp),\bu(t), t, \bp)\\[0.5em]
&\be^L \le \be\big(\bx(t_0, \bp), \bx(t_f, \bp), t_0, t_f, \bp\big)  \le \be^U \\[0.5em]
&\bh^L \le  \bh(\bx(t, \bp),\bu(t), t, \bp)  \le \bh^U \\
\end{array}
\right. &  \nonumber 
\end{eqnarray}
Problem~$(A^\infty)$ generates a large family of transcendental tychastic optimal control problems parameterized by the cost functional $J_{tyc}$.  As noted before, there is an absence of a generic tychastic cost functional. A nonexhaustive collection of possible tychastic functionals $J_{tyc}$ are given by \eqref{eq:minimax}, \eqref{eq:tyc-avg}, \eqref{eq:tyc-trCov}, \eqref{eq:tyc-E-avg} and \eqref{eq:tyc-Pr}.

\begin{figure}[!ht]
	\centering
    \includegraphics[angle = 0, width=0.75\columnwidth,clip]{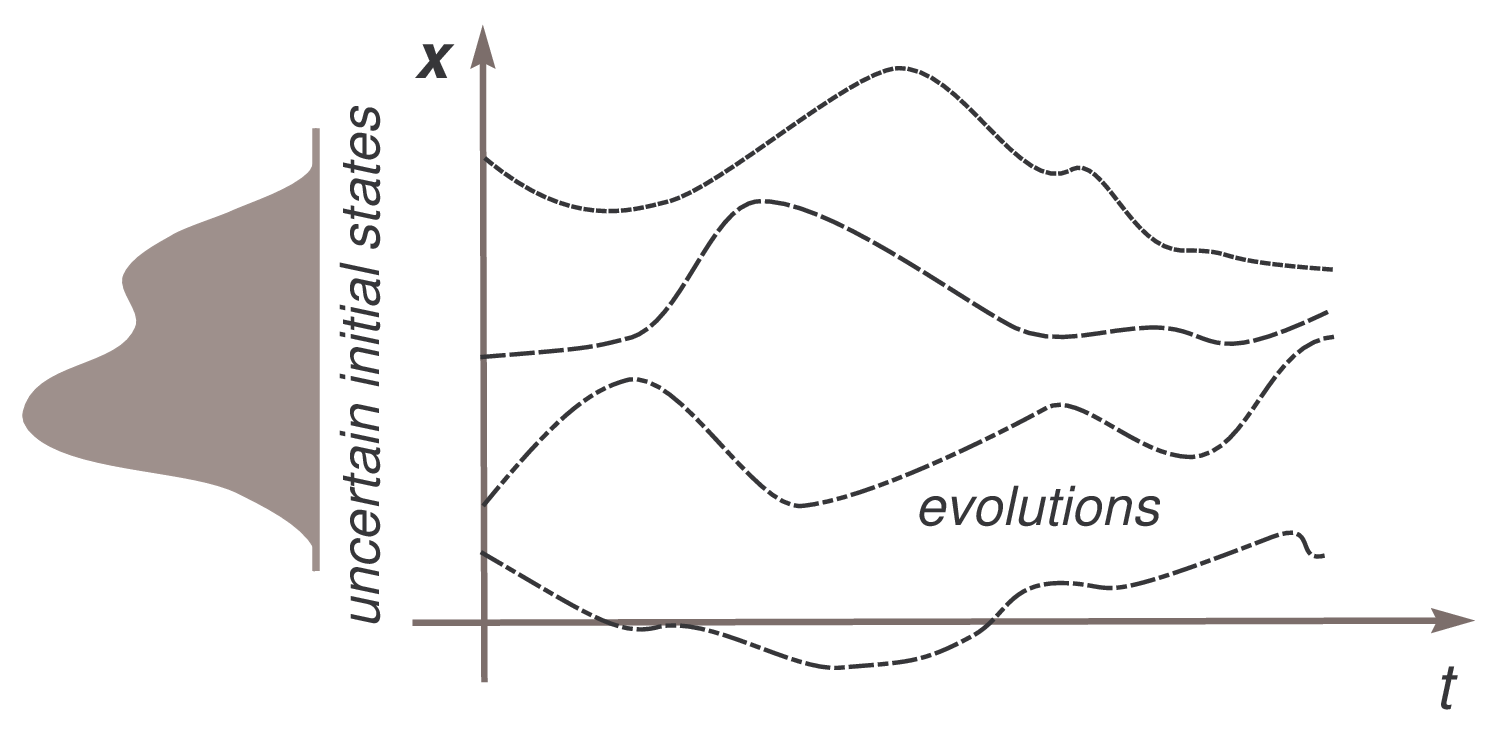}
	\caption{Schematic for the evolutions of state trajectories for uncertain initial conditions; figure adapted from [\citen{LS:acc-2}].}
	\label{fig:RandomIC}
\end{figure}

Note that \emph{\textbf{the random variable $\bp$ can appear anywhere in a tychastic trajectory optimization problem}}; for instance, the initial state may be random, in which case, we set $N_x$-components of $\bp$, say $\bp^0$ to be $\bx(t_0)$; i.e., $\bp^0 = \bx(t_0)$.  In this instance, the schematic of Fig.~\ref{fig:DiffInc} is modified as shown in Fig.~\ref{fig:RandomIC}.
%

The severity of the non-dynamical constraints in Problem~$(A^\infty)$ is that it requires the entire ensemble of solutions (generated by each $\bp \in supp(\bp)$) to satisfy all constraints.  Although this is theoretically possible under certain situations\cite{uoc-2:SD,ross-book}, this fact is not known a priori.  Hence, it is necessary to construct the most desirable solution by way of constructing a sequence of ``simpler'' problems whose solutions inform the construction of the next problem in the sequence.  Conceptually, each ancillary problem in this sequence approaches Problem~$(A^\infty)$.  To put this in perspective, consider a particular version of a chance-constrained\cite{charnes-cooper} tychastic optimal control problem framed as\cite{LS:acc-2,ross-book},
\begin{eqnarray}
&(C) \left\{
\begin{array}{lll}
\text{Minimize }    J_{tyc}[(\bx(\cdot, \cdot), t_0, t_f; supp(\bp))] \\[0.5em]
\text{Subject to }   \\[0.5em]
 \dot{\bx}(t, \bp) = {\bff}(\bx(t, \bp),\bu(t), t, \bp) \quad \forall\ \bp \in supp(\bp)\\[0.5em]
\text{Pr}\set{\be^L \le \be\big(\bx(t_0, \bp), \bx(t_f, \bp), t_0, t_f, \bp\big)  \le \be^U} \ge 1 - r_e \\[0.5em]
\displaystyle\max_{t \in [t_0, t_f]}\text{Pr}\set{\bh^L \le  \bh(\bx(t, \bp),\bu(t), t, \bp)  \le \bh^U} \ge 1 - r_h\\
\end{array}
\right. &  \label{eq:prob-C}
\end{eqnarray}
where, $r_e \in (0,1) $ is the allowable risk in not satisfying any of the boundary conditions and $r_h \in (0,1)$ is the maximum allowable risk in not satisfying any of the path constraints. Ignoring the fact that even computing the constraints in Problem~$(C)$ is not easy (see \eqref{eq:Pr-def}) it is not clear how to specify the risk levels.  For instance, if the risk level is stipulated to be close to zero, then it is apparent that a solution to Problem~$(C)$ might not exist because of the demand of an unachievable low value of risk.

In order to construct the most desirable and feasible tychastic trajectory, we construct a sequence of tractable tychastic optimal control problems wherein each element of the sequence is based on the solution to the prior problem in the sequence\cite{ross-book}.  Thus, problem formulation and solution generation go hand in hand.  An unscented solution to a tychastic optimal control problem forms the critical computational tool for constructing this sequence of problems and solutions.

\section{Unscented Techniques for Tychastic Trajectory Optimization}
\label{sec:tyc2unscent}
Let $(\bp_i, w_i)$ be a collection of cubature ``nodes'' and ``weights'' that approximate a Lebesgue-Stieltjes integral according to\cite{engels,cools_2001},
\begin{equation}\label{eq:cubature}
\int\cdots\int_{supp(\bp)} \bg(\bp)\, dm(\bp) = \lim_{n \to \infty} \sum_{i=1}^n w_i \bg(\bp_i)
\end{equation}
where, $\bg: \bp \mapsto \real{N_g}$ is an integrable function.  Then for any finite $n$, we define a \emph{semi-discretized tychastic optimal control problem} by replacing all the expectation operators by their cubature schemes.  Thus, for example, for deterministic clock times, Problem~$(A^\infty)$ generates the following semi-discretized, large-scale, deterministic optimal control problem,
\begin{eqnarray}
& \forall\ \bp_i, \ i = 1, \ldots, n & \nonumber\\[0.5em]
&(A^n) \left\{
\begin{array}{lll}
\text{Minimize }   &J^n_{tyc}[(\bX(\cdot), \buf, t_0, t_f; \bP)]  \\[0.5em]
\text{Subject to }  & \dot{\bX}(t) = {\bff}(\bX(t),\bu(t), t, \bP)\\[0.5em]
&\be^L \le \be\big(\bx(t_0, \bp_i), \bx(t_f, \bp_i), t_0, t_f, \bp_i\big)  \le \be^U \\[0.5em]
&\bh^L \le  \bh(\bx(t, \bp_i),\bu(t), t, \bp_i)  \le \bh^U \\
\end{array}
\right. & \nonumber 
\end{eqnarray}
where, $\bX(t) \in \real{n\,N_x}$  and $\bP \in \real{N_p \times n}$ are given by
\begin{equation}
\bX(t) := (\bx(t, \bp_1), \ldots, \bx(t, \bp_n)), \quad \bP := (\bp_1, \ldots, \bp_n)
\end{equation}
In the definition of Problem~$(A^n)$, the function $J^n_{tyc}$ is assumed to be computed in a manner consistent with the statistics that is being minimized in Problem~$(A^\infty)$. See, for example, \eqref{eq:minimax}, \eqref{eq:tyc-avg}, \eqref{eq:tyc-trCov}, \eqref{eq:tyc-E-avg} and \eqref{eq:tyc-Pr}. Furthermore, we have also assumed that the function $\bff$ is overloaded according to,
\begin{equation}
\bff(\bX(t),\bu(t), t, \bP) := \big(\bff(\bx(t, \bp_1),\bu(t), t, \bp_1), \ldots, 
\bff(\bx(t, \bp_n),\bu(t), t, \bp_n)\big)
\end{equation}
In unscented trajectory optimization, we simply use the family of $n = N_\sigma$ sigma points of Julier et al\cite{julier-acc-95,julier:simplex} or their various generalizations\cite{UT-various-2012,UT-various-2015,UT-various-2021,uo-1:MRYpts-2016}.  Because $N_\sigma$ is substantially smaller than the number of Monte Carlo points, the resulting scale of the semi-discretized tychastic optimal control problem is correspondingly smaller.  Similar small scale unscented trajectory optimization problems may be constructed based on their tychastic counterparts.

The sigma points may be viewed as low-order cubature points\cite{cools_2001,UT-various-2021} that match low-order moments; hence, unscented trajectory optimization techniques do not adequately capture higher-order statistics such as the tail of the distribution.  Consequently, if an unscented trajectory optimization problem is solved, its statistics must be verified or recomputed using Monte Carlo simulations, in addition to performing the usual verification and validation of the computed optimal solutions\cite{ross-book}.  \emph{\textbf{That is, Monte Carlo methods are used/useful for simulation while unscented methods are used/useful for trajectory optimization.}}

Coupled with the fact that an unscented trajectory optimization problem that corresponds to a particular tychastic optimal control problem may not have a solution, a systematic procedure is needed to determine what exactly is statistically controllable and by what margins.  Thus, for example, although a chance-constrained optimal control problem (cf.~\eqref{eq:prob-C}) is part of the collection of tychastic optimal control problems, a specification of a small value of risk may not be feasible whereas a slightly higher one might make multiple trajectory solutions possible.  Such mission design questions can be quickly answered by solving an appropriate sequence of unscented trajectory optimization problems.  The particular sequence of problems that needs to be solved depends upon the mission requirements and the origins of the uncertainties.

A general procedure to formulate and solve the the most desirable tychastic trajectory optimization problem is posed as follows:
\begin{enumerate}
\item[\textbf{Step 1}] Formulate a deterministic trajectory optimization problem using the most probable value of the uncertain parameters, $\bp = \bp^b$.  This generates a ``baseline'' deterministic optimal control problem:
\begin{eqnarray}
&(B) \left\{
\begin{array}{lll}
\text{Minimize }   &J^{B}[\bx(\cdot), \buf, t_0, t_f] \\[0.5em]
\text{Subject to }  & \dot{\bx}(t, \bp^b) = {\bff}(\bx(t, \bp^b),\bu(t), t, \bp^b)\\[0.5em]
&\be^L \le \be\big(\bx(t_0, \bp^b), \bx(t_f, \bp^b), t_0, t_f, \bp^b\big)  \le \be^U \\[0.5em]
&\bh^L \le  \bh(\bx(t, \bp^b),\bu(t), t, \bp^b)  \le \bh^U \\
\end{array}
\right. &  \nonumber
\end{eqnarray}
Solve Problem~$(B)$ (by any method). In the examples to follow we used DIDO$^\copyright$\cite{DIDO:arXiv}, a MATLAB$^\circledR$ optimal control toolbox\footnote{\textit{https://www.mathworks.com/products/connections/product\_detail/dido.html.  See also https://elissarglobal.com/}} that implements the spectral algorithm\cite{spec-alg,auto-knots,RossKarp_IFAC_2012} for pseudospectral (PS) methods\cite{RossKarp_IFAC_2012,arb-grid,fastmesh,BirkhoffTN}.
Verify and validate (V\&V) this result for feasibility and optimality using standard procedures\cite{ross-book,DIDO:arXiv}. Note that the verification process is critical for the remainder of the steps because, for instance, if the deterministic solution is not independently feasible\cite{ross-book}, then the Monte Carlo analysis in Step~\textbf{2.} will generate incorrect results. Let $t \mapsto \bu^b(t)$ be the optimized control trajectory obtained in this step.
\item[\textbf{Step 2}] Using a Runge-Kutta solver\cite{hnw-ode} and the tychastic ODE $\dot\bx = \bff(\bx, \bu^b(t), t, \bp)$ generated by $\bu^b(\cdot)$, initiate a Monte Carlo simulation over the uncertain parameters, $\bp \in supp(\bp)$.  Analyze the results for various statistics such as violations of the endpoint constraints and path constraints. If the results are satisfactory, stop.  The problem is insensitive to the uncertain parameters.
\item[\textbf{Step 3}] Based on the results \textbf{Step~2}, develop a preliminary tychastic trajectory optimization that addresses managing the statistics of some element of the constraint violations.  For example, if \textbf{Step~2} indicates that the average values of the endpoint constraints and/or path constraints are violated then it may be necessary to formulate the following ancillary tychastic problem\cite{ug:acc,ross-book}:
\begin{eqnarray}
& \forall\ \bp \in supp(\bp) & \nonumber\\[0.5em]
&(A_1) \left\{
\begin{array}{lll}
\text{Minimize }   &J_{tyc}[(\bx(\cdot, \cdot), \buf, t_0, t_f, supp(\bp))]  \\[0.5em]
\text{Subject to }  & \dot{\bx}(t, \bp) = {\bff}(\bx(t, \bp),\bu(t), t, \bp)\\[0.5em]
&\be^L \le \mathcal{E}\big[\be\big(\bx(t_0, \bp), \bx(t_f, \bp), t_0, t_f, \bp\big)\big]  \le \be^U \\[0.5em]
&\bh^L \le  \mathcal{E}\big[\bh(\bx(t, \bp),\bu(t), t, \bp) \big]  \le \bh^U \\
\end{array}
\right. & \nonumber
\end{eqnarray}
Using \eqref{eq:cubature}, the unscented instantiation of Problem~$(A_1)$ is given by,
\begin{eqnarray}
& \forall\ \bp_i, \ i = 1, \ldots, n & \nonumber\\[0.5em]
&(A_1^U) \left\{
\begin{array}{lll}
\text{Minimize }   &J^n_{tyc}[(\bX(\cdot), \buf, t_0, t_f, \bP)]  \\[0.5em]
\text{Subject to }  & \dot{\bX}(t) = {\bff}(\bX(t),\bu(t), t, \bP)\\[0.5em]
&\be^L \le \sum_{i=1}^{N_\sigma} w_i\be\big(\bx(t_0, \bp_i), \bx(t_f, \bp_i), t_0, t_f, \bp_i\big)  \le \be^U \\[0.5em]
&\bh^L \le  \sum_{i=1}^{N_\sigma} w_i\bh(\bx(t, \bp_i),\bu(t), t, \bp_i)  \le \bh^U \\
\end{array}
\right. & \nonumber
\end{eqnarray}
\item[\textbf{Step 4}] Solve the unscented trajectory optimization problem developed in \textbf{Step~3}. Because this unscented problem is deterministic, it should be solvable by the method used in \textbf{Step~1}. Let $t \mapsto \bu^1(t)$ be the optimized unscented control trajectory obtained in this step. V\&V this result as performed in \textbf{Step~1}.
\item[\textbf{Step 5}] Using the ODE $\dot\bx = \bff(\bx, \bu^1(t), t, \bp)$, conduct a Monte Carlo analysis over the uncertain parameters, $\bp \in supp(\bp)$. If the results are satisfactory, stop; else go to \textbf{Step~6}.
\item[\textbf{Step 6}] Based on the results of \textbf{Step~5}, develop a second ancillary tychastic trajectory optimization that further controls the statistics of one or more element of the constraint violations.  For example, if \textbf{Step~5} indicated that the average values of the endpoint constraints and/or path constraints were successfully controlled, it may be more desirable to shrink their dispersions. Suppose that it was necessary to reduce the dispersions of $\bx(t_f, \bp)$.  Then, the second ancillary tychastic trajectory optimization problem may be posed using \eqref{eq:tyc-trCov}.  From \eqref{eq:cubature}, \eqref{eq:expect-def} and \eqref{eq:cov-def}, the unscented instantiation of \eqref{eq:tyc-trCov} can be written as,
 \begin{equation}
 \text{tr} \left[ \sum_{i=1}^{N_\sigma}w_i\left[\bx(t_f, \bp_i) - \overline{\bx}_f \right] \left[\bx(t_f, \bp_i) - \overline{\bx}_f \right]^T \right]
 \end{equation}
where,
 \begin{equation}
 \overline{\bx}_f :=  \sum_{i=1}^{N_\sigma}w_i\,\bx(t_f, \bp_i)
 \end{equation}
\item[\textbf{Step 7}] Solve the unscented trajectory optimization problem posed in \textbf{Step~6} and repeat \textbf{Steps~4}, \textbf{5} and \textbf{6} until a desirable solution is found.
\end{enumerate}

We emphasize at this point that the tychastic optimal control problems defined in the preceding steps are not the only possible ones.  In fact, a large number of additional problems may be formulated based on the specific needs of a mission; see for example \cite{LS:acc-1,LS:acc-2,ross-book,uo-1:2015,uoc-2:SD,uoc-patent-1,uoc-patent-2} where many number of other tychastic cost functionals and constraints are presented based on the definition of success.

Note also that \textbf{Steps 1}, \textbf{4} and \textbf{7} involve solving a deterministic optimal control problem. Because the theoretical analysis of the existence of a solution to deterministic optimal control problem is not quite straightforward\cite{vinter}, it is necessary to employ certain alternative methods to determine if solutions exist and are valid.  Most of these methods are based on V\&V techniques that involve feasibility and optimality tests which are generated by an application of Pontryagin's Principle\cite{ross-book}.  See \cite{RossKarp_IFAC_2012,ross-book,DIDO:arXiv,IEEE:spectrum,Bedrossian_2009,SIAMnews,LRO-CSM} for an extensive discussion of these tests, most notably as they apply to real-world flight implementations\cite{RossKarp_IFAC_2012,IEEE:spectrum,Bedrossian_2009,SIAMnews, LRO-CSM}.  Some of these techniques are illustrated in the sections to follow.

\section{Tychastic Trajectory Optimization Illustrated Via the Unscen\-ted Transform }
In this section, we illustrate the ideas presented in the prior sections using various formulations of the Zermelo problem\cite{serres}. The classic deterministic Zermelo problem is discussed in many optimal control textbooks\cite{bryson:ho,longuski,ross-book} because of its richness in illustrating optimal control concepts\cite{trelat:survey}. Consequently, we formulate different versions of the Zermelo problems to offer similar insights.

\subsection{Step 1: A Baseline Deterministic Zermelo Problem Formulation and Its Solution}
\begin{figure}[!h]
      \begin{center}
        \includegraphics[angle = 0,width = 0.48\columnwidth, clip] {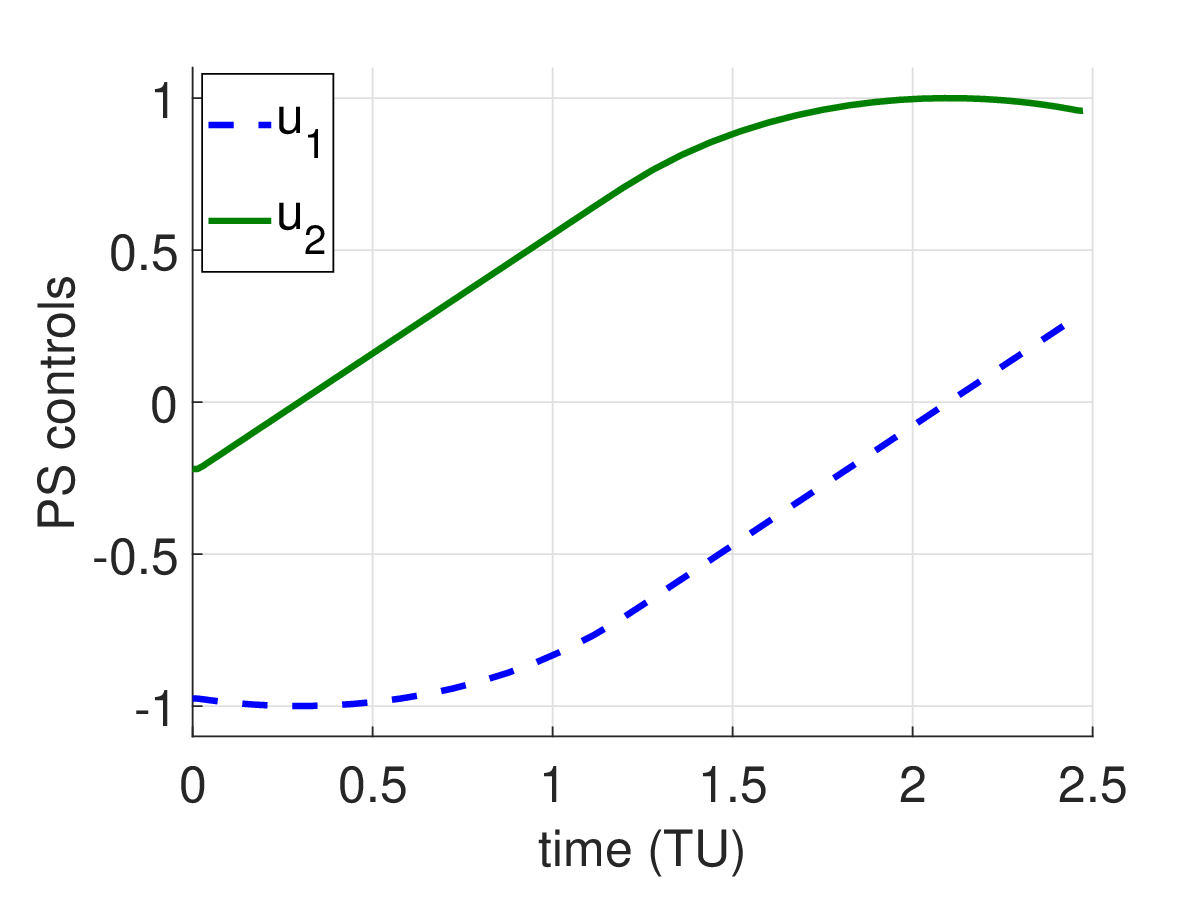}
        \includegraphics[angle = 0,width = 0.48\columnwidth, clip] {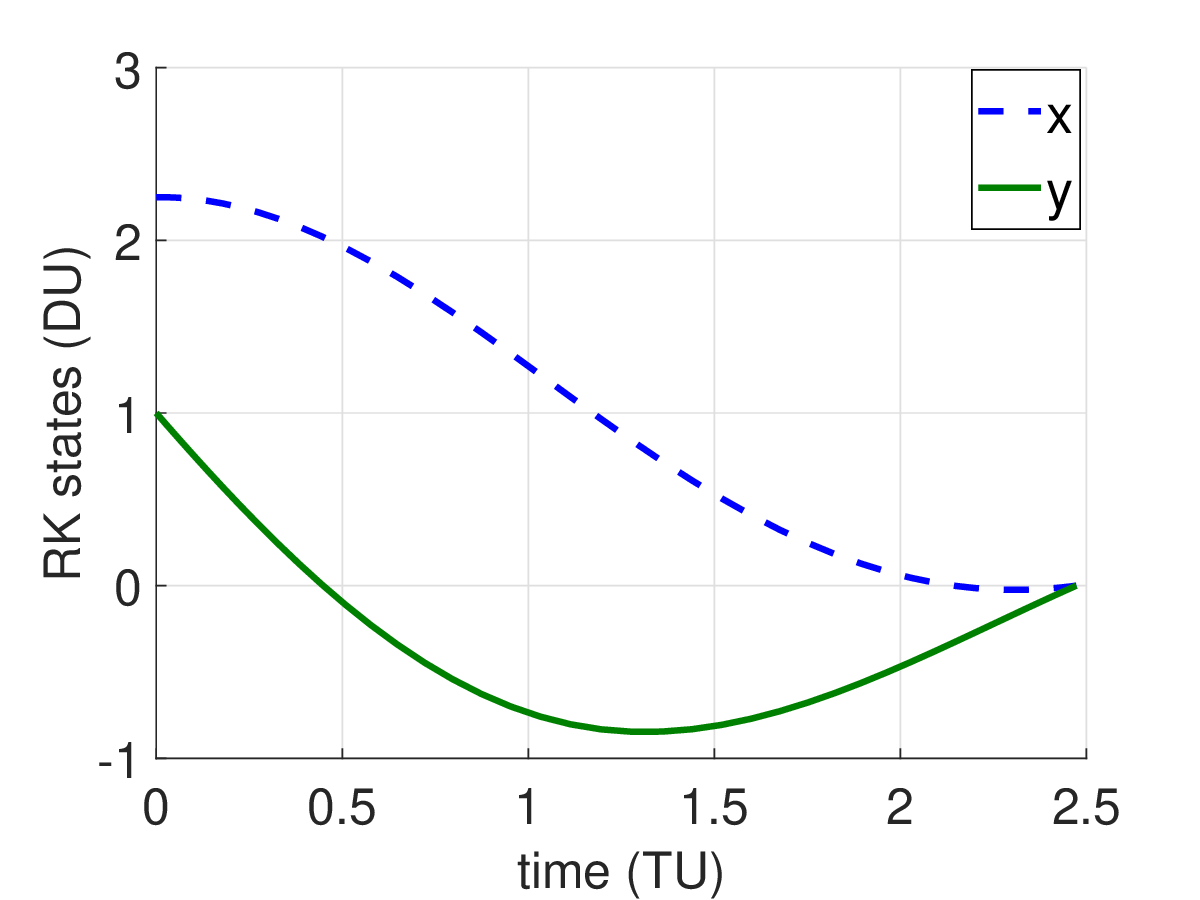}
      \caption{A candidate optimal solution to Problem $Z_0$.}
      \label{fig:Zsol0}
      \end{center}
\end{figure}

According to the historical and technical review by Serres\cite{serres}, E. Zermelo formulated his famous problem in 1931 in connection with the steering of a ship in a wind vector field. Since then, the problem has been connected to other problems in disparate fields including many in aerospace engineering: North-South stationkeeping of a spacecraft\cite{kamel,kech}, control of an electrodynamic tether\cite{stevens-Zerm}, and proximity operations for relative motion\cite{uoc-2:SD}.  The differential equations for Zermelo's problem can be written as,
\begin{equation}
\begin{aligned}
\dot x &= W_1(x, y, t) + u_1 \\
\dot y &= W_2(x, y, t) + u_2
\end{aligned} \label{eq:Z-1}
\end{equation}
where $\mathbf{W}(x,y) := (W_1(x, y), W_2(x,y))$ is the wind vector field, and $\bu = (u_1, u_2)$ is the steering vector normalized by the ship's speed and constrained by $u_1^2 + u_2^2 = 1$.  Choosing a $(p, q)$-parameterized linear cross-wind profile, $\mathbf{W}(x,y) := (py, qx)$, a deterministic Zermelo problem with baseline values of $(p^b, q^b) = (1, -1)$ is posed as,
\begin{eqnarray*}
&\bx := (x, y) \in  \real{2} & \nonumber\\
& \bu \in \U := \set{(u_1, u_2) \in \real{2}:\ u_1^2 + u^2_2 = 1 } & \nonumber\\
& \bp^b := (p^b, q^b) = (1, -1) \in \real{2}  & \nonumber\\
&(Z_0) \left\{
\begin{array}{lll}
\text{Minimize }  & J^b[\bx(\cdot), \buf, t_f]
:= t_f  \\
\text{Subject to} & \dot x(t) =  p^b y(t) + u_1(t)\\
&\dot y(t) = q^b x(t) + u_2(t)\\
&(\bx(t_0), t_0) = (2.25, 1, 0) \\
&\bx(t_f) =  (0,0) \\
\end{array} \right.& 
\end{eqnarray*}
A baseline optimal solution to this problem (i.e. control-state pair: $t \mapsto (\bu^b(t), \bx^b(t))$) is shown in Fig.~\ref{fig:Zsol0}.
%
%
As noted in Section~\ref{sec:tyc2unscent}, candidate optimal solutions to all trajectory optimization problems were obtained using DIDO\cite{DIDO:arXiv,ross-book}. These DIDO-generated controls are labeled PS controls in Fig.~\ref{fig:Zsol0}. In accordance with \textbf{Step~1} of the procedure outlined in Section~\ref{sec:tyc2unscent}, the state trajectory shown in Fig.~\ref{fig:Zsol0} was generated by a Runge-Kutta propagation (RK45 method implemented in \texttt{\textsf{ode45}} in MATLAB) of the initial condition using linearly interpolated values of the PS controls.  Thus, the state trajectory generation is agnostic to the DIDO/PS controls and constitutes part of an independent V\&V. This two-step process of generating verifiable solutions is extensively discussed elsewhere\cite{RossKarp_IFAC_2012,ross-book,DIDO:arXiv} and is the same standardized procedure used in the pre-flight validation of DIDO solutions implemented in many NASA missions\cite{RossKarp_IFAC_2012,IEEE:spectrum,Bedrossian_2009,SIAMnews,LRO-CSM}. Further details of this two-step procedure for generating flight-implementable optimal solutions are described in [\citen{DIDO:arXiv}] and [\citen{RossKarp_IFAC_2012}].

It is apparent from Fig.~\ref{fig:Zsol0} that the PS controls are separately/independently verified for feasibility because the RK45-propagated trajectories indeed acquire the targeted value of $(x_f, y_f) = (0,0)$ (within specified RK tolerances).  Optimality verification tests were also performed on this result in accordance with the methods discussed in [\citen{ross-book}] and [\citen{DIDO:arXiv}]; they are not discussed here for the purposes of brevity.  The baseline value of the cost functional (i.e., minimum transfer time) was found to be $t^b_f \approx 2.47$.

\subsection{Step 2: Monte Carlo Analysis for the Baseline Solution}
\begin{figure}[!h]
      \begin{center}
        \includegraphics[angle = 0,width = 0.48\columnwidth] {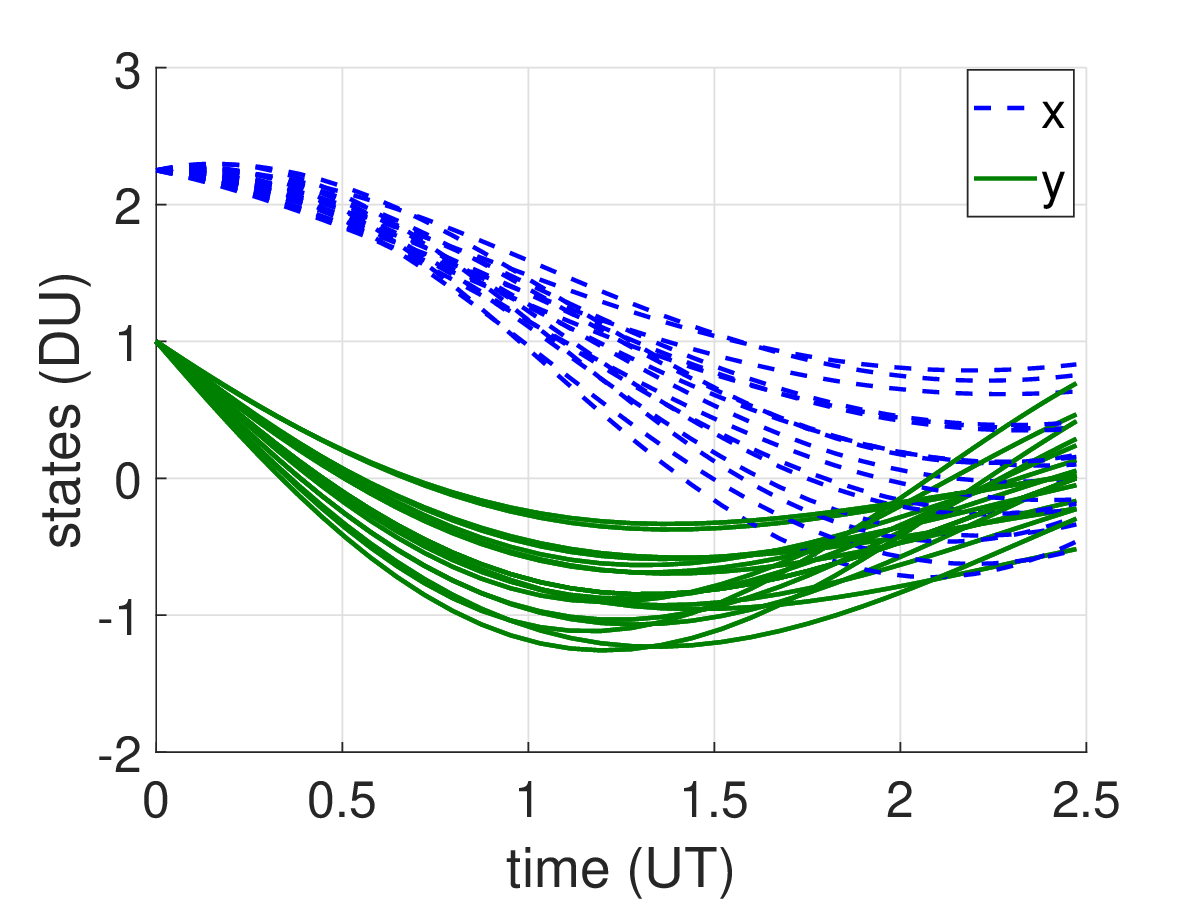}
        \includegraphics[angle = 0,width = 0.48\columnwidth] {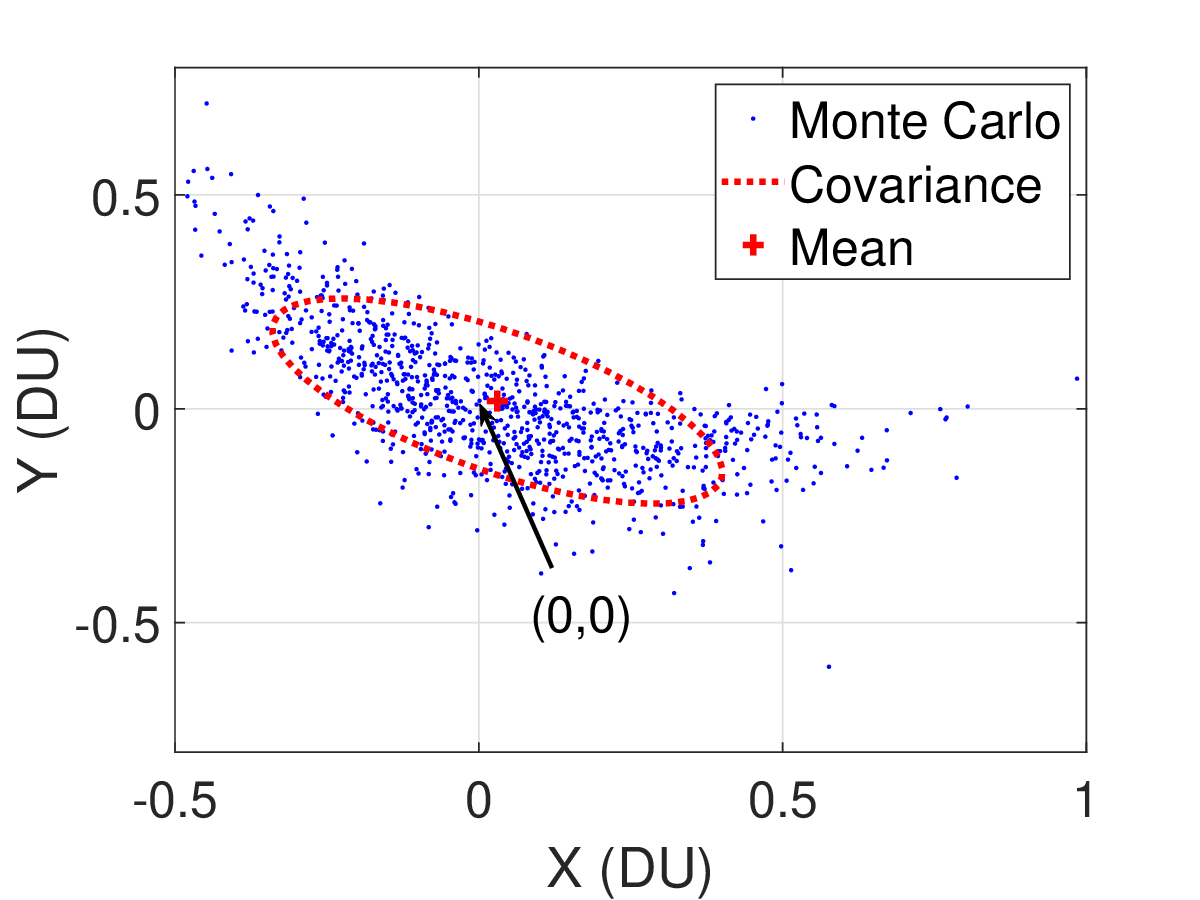}
      \caption{Results of a Monte-Carlo simulation for $\bp$ given by \eqref{eq:p=PDFdata}.  See also Fig.~\ref{fig:DiffInc} and compare with Fig.~\ref{fig:stocTrajs}.}
      \label{fig:MC4Z0}
      \end{center}
\end{figure}

To perform \textbf{Step~2} of the procedure outlined in Section~\ref{sec:tyc2unscent}, we now suppose that $\bp$ is a random parameter with a Gaussian distribution given by,
\begin{equation}\label{eq:p=PDFdata}
\bp:= (p, q) \sim \mathcal{N}\big((1, -1), diag(0.2^2, 0.1^2)\big)
\end{equation}
A 1000-point Monte Carlo simulation for this uncertain data is shown in Fig.~\ref{fig:MC4Z0}.
%
%
For clarity of presentation, the state trajectories shown in Fig.~\ref{fig:MC4Z0} are for a smaller sample set (than 1000 simulations).  In any event, it is apparent that nearly all of the endpoints do not meet the target condition of $\bx_f = (0,0)$.  In other words, we have,
\begin{equation}\label{eq:Pr=040}
\text{Pr}\set{\bx(t_f, \bp) = (0,0)} \approx 0
\end{equation}
Even the mean of the final states (as computed by using the results of the Monte Carlo simulations) is not at the origin; i.e.,
\begin{equation}\label{eq:meanZermNE0}
\int_{-\infty}^{\infty}\int_{-\infty}^{\infty} \bx(t_f, \bp)\, d\Phi\left(\frac{p-1}{0.2}\right) d\Phi\left(\frac{q+1}{0.1}\right) \approx (0.03, 0.02)  \neq (0, 0)
\end{equation}
where,
$$ \Phi(\xi) :=  \frac{1}{\sqrt{2\pi}} \int_{-\infty}^\xi e^{-s^2/2}ds $$

\subsection{Step 3: An Ancillary Unscented Zermelo Problem Formulation}
Because the mean value of the final conditions is not at the origin (cf.~\eqref{eq:meanZermNE0}), consider now an ancillary tychastic optimal control problem which requires the mean value of the final time conditions to be at the origin.  This simple requirement translates to the following problem:
\begin{eqnarray*}
&\bu \in \U := \set{(u_1, u_2) \in \real{2}:\ u_1^2 + u^2_2 = 1 } & \nonumber\\
& \bp:= (p, q) \sim \mathcal{N}\big(1, -1), diag(0.2^2, 0.1^2)\big)   & \nonumber\\
&(Z_1) \left\{
\begin{array}{lll}
\text{Minimize }  & J^1_{tyc}[\bx(\cdot, \cdot), \buf, t_f; supp(\bp)]
:= t_f  \\
\text{Subject to} & \dot x(t,\bp) =  p y(t, \bp) + u_1(t) \\
&\dot y(t, \bp) = q x(t, \bp) + u_2(t)\\
&(\bx(t_0, \bp), t_0) = (2.25, 1, 0) \\[0.5em]
&\displaystyle\int_{-\infty}^{\infty}\int_{-\infty}^{\infty} x(t_f, \bp) d\Phi\left(\frac{p-1}{0.2}\right) d\Phi\left(\frac{q+1}{0.1}\right)= 0 \\[1em]
&\displaystyle\int_{-\infty}^{\infty}\int_{-\infty}^{\infty} y(t_f, \bp) d\Phi\left(\frac{p-1}{0.2}\right) d\Phi\left(\frac{q+1}{0.1}\right) = 0 \\
\end{array} \right.& 
\end{eqnarray*}
Let,
$\big((p_1, q_1), \ldots, (p_{N_\sigma}, q_{N_\sigma})\big)$ be the sigma points\cite{julier:simplex,julier-acc-95} or its generalizations\cite{UT-various-2012,UT-various-2015,UT-various-2021,uo-1:MRYpts-2016} associated with $\mathcal{N}\left((1, -1)\right.$, $\left.diag\left(0.2^2, 0.1^2\right)\right)$.  Then, the unscented instantiation of Problem~$(Z_1)$ is given by,
\begin{eqnarray*}
& \bu \in \U := \set{(u_1, u_2) \in \real{2}:\ u_1^2 + u^2_2 = 1 } & \nonumber\\
& \bX := \set{(x_i, y_i), i = 1, \ldots, N_\sigma} \in \real{2N_\sigma}   & \nonumber\\
&(Z_1^U) \left\{
\begin{array}{lll}
\text{Minimize }  & J^1_{U}[\bX(\cdot), \buf, t_f]
:= t_f  \\
\text{Subject to} & \dot x_i(t) =  p_i y(t) + u_1(t) &\forall\ i =1, \ldots, N_\sigma\\
&\dot y_i(t) = q_i x(t) + u_2(t)  &\forall\ i =1, \ldots, N_\sigma\\
&(x_i(t_0), y_i(t_0)) = (2.25, 1)   &\forall\ i =1, \ldots, N_\sigma \\
&t_0 = 0 \\
&\sum_{i=1}^{N_\sigma} w_i\, x_i(t_f) = 0 \\
&\sum_{i=1}^{N_\sigma} w_i\, y_i(t_f) = 0 \\
\end{array} \right.& 
\end{eqnarray*}
where $w_i, \ i = 1, \ldots, N_\sigma$ are the weights associated with the selected sigma points.  Because Problem~$(Z_1^U)$ is deterministic, it can be solved using DIDO as before.

\subsection{Steps 4 and 5: Solve and Analyze Problem~$(Z_1^U)$}
\begin{figure}[!h]
      \begin{center}
        \includegraphics[angle = 0,width = 0.48\columnwidth] {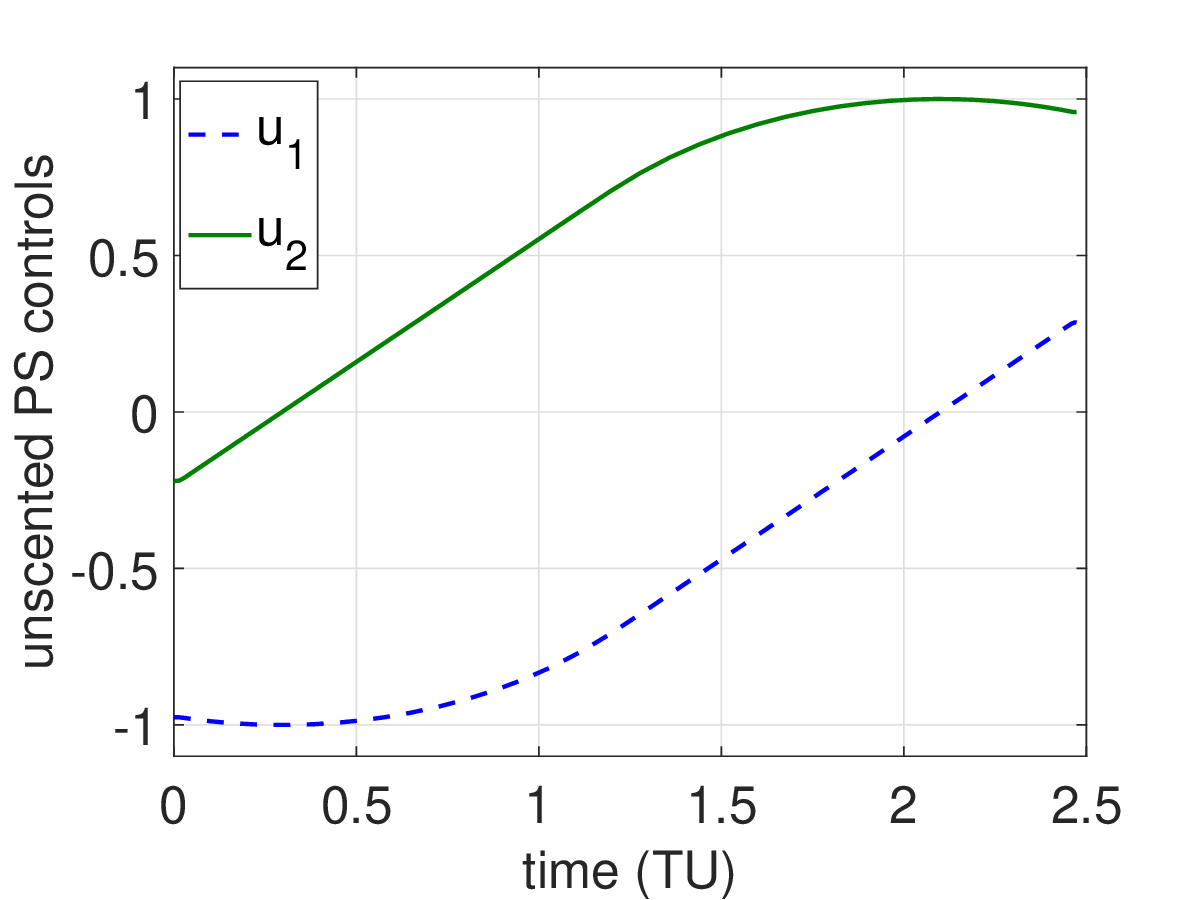}
        \includegraphics[angle = 0,width = 0.48\columnwidth] {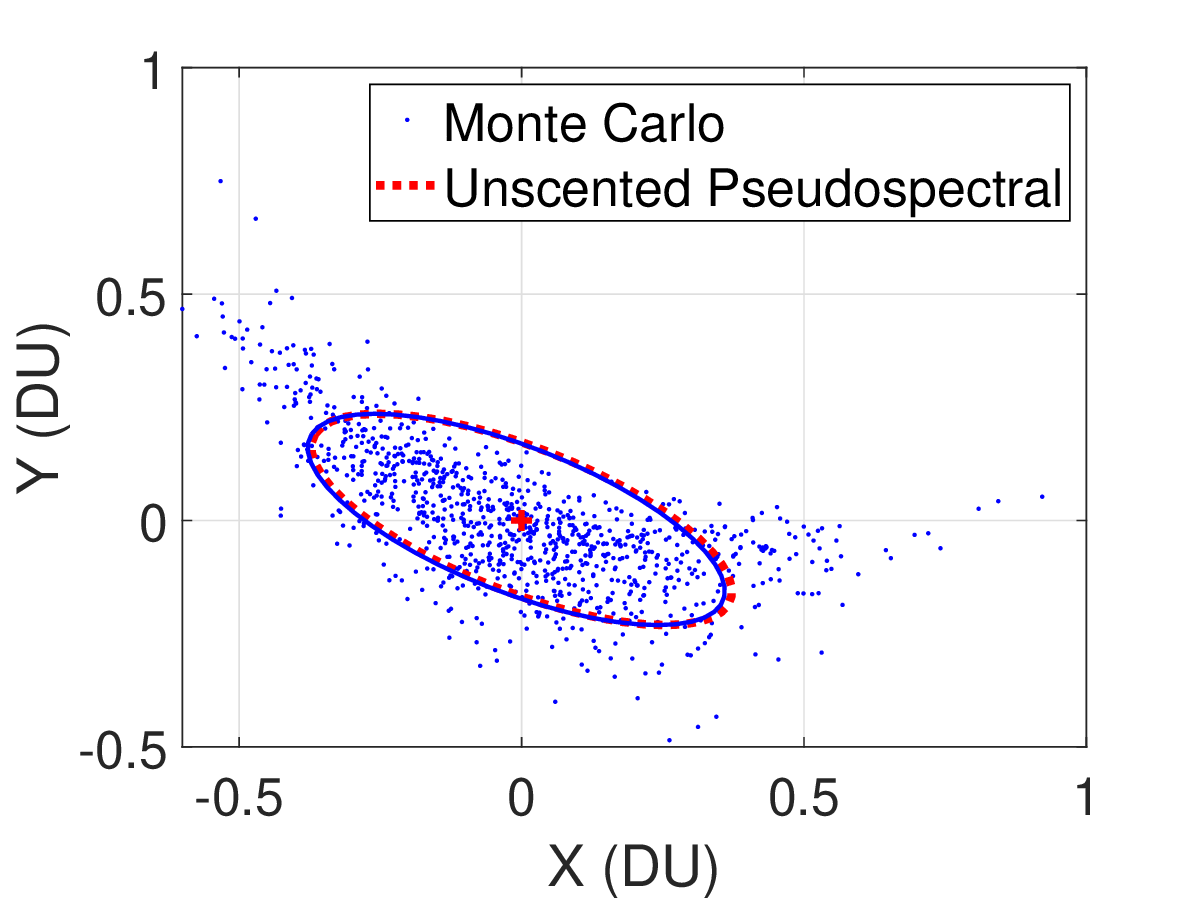}
      \caption{An unscented PS solution to Problem $Z_1$ and its Monte Carlo results.}
      \label{fig:Z1Sol}
      \end{center}
\end{figure}

A PS control solution to Problem~$(Z_1^U)$ is shown in Fig.~\ref{fig:Z1Sol}.
%
%
These unscented controls were obtained using DIDO in the same manner as in \textbf{Step~1} and verified accordingly. By visual inspection, it is apparent that the unscented PS control solution looks nearly the same as the deterministic one (compare Fig.~\ref{fig:Zsol0}). Furthermore the optimal cost (i.e., transfer time) for Problem~$(Z_1^U)$ was found to be about $2.7$.  Comparing this value to the result from \textbf{Step~1}, we have,
\begin{equation}
J^b[\bx^b(\cdot), \bu^b(\cdot), t^b_{f}] \approx J^1_{tyc}[\bx^1(\cdot, \cdot), \bu^1(\cdot), t^1_{f}; supp(\bp)]  \approx 2.47
\end{equation}
Hence, one can argue that targeting the mean value to the origin is relatively free of charge.

Next, in accordance with \textbf{Step~5}, Monte Carlo simulations were performed using the ODE $\dot\bx = \bff(\bx, \bu^1(t), t, \bp)$, where, $\bp \sim \mathcal{N}\big(1, -1), diag(0.2^2, 0.1^2)\big)$ and $\bu^1(\cdot)$ is the unscented control trajectory shown in Fig.~\ref{fig:Z1Sol}.  The results of this Monte Carlo simulation are shown in Fig.~\ref{fig:Z1Sol}. Although it is not visually apparent, the Monte Carlo simulations verify that the mean value has indeed shifted to the origin as required. It was numerically verified to be $(0,0)$.  These tests indicate that the tychastic trajectory optimization problem given by Problem~($Z_1$) was successfully solved using its unscented instantiation given by Problem~$(Z_1^U)$.

\subsection{Step 6: Controlling the Dispersions of the Final Positions}
Having successfully driven the mean values of the states to the origin in \textbf{Step~5}, consider now levying an additional requirement of reducing the dispersion of the endpoint values.  There are several challenges in translating this engineering requirement to a precise mathematical problem formulation.  The first challenge is simply a precise definition of dispersion. The second challenge involves a real number (a measure) that defines a tolerable value of dispersion.  As noted before, it is intuitively obvious that if the tolerance is set too tight, a solution might not exist.  As a means to address this challenge, we simply change the tychastic cost functional to minimize dispersions. Furthermore, we simply take the trace of the covariance matrix as a measure of dispersion.  These stipulations lead us to construct the following tychastic optimal control problem:
\begin{eqnarray*}
&\bx := (x, y) \in  \real{2} \qquad \bu \in \U := \set{(u_1, u_2) \in \real{2}:\ u_1^2 + u^2_2 = 1 } & \nonumber\\
& \bp:= (p, q) \sim \mathcal{N}\big(1, -1), diag(0.2^2, 0.1^2)\big)   & \nonumber\\
&(Z_2) \left\{
\begin{array}{lll}
\text{Minimize }  & J^2_{tyc}[\bx(\cdot, \cdot), \buf, t_f; supp(\bp)] := \text{tr}\,\text{Cov}[\bx(t_f, \bp)]  \\
\text{Subject to} & \dot x(t, \bp) =  p\, y(t, \bp) + u_1(t)\\
&\dot y(t, \bp) = q\, x(t, \bp) + u_2(t)\\
&(\bx(t_0, \bp), t_0) = (2.25, 1, 0) \\[0.5em]
&\displaystyle\int_{-\infty}^{\infty}\int_{-\infty}^{\infty} x(t_f, \bp) d\Phi\left(\frac{p-1}{0.2}\right) d\Phi\left(\frac{q+1}{0.1}\right)= 0 \\[1em]
&\displaystyle\int_{-\infty}^{\infty}\int_{-\infty}^{\infty} y(t_f, \bp) d\Phi\left(\frac{p-1}{0.2}\right) d\Phi\left(\frac{q+1}{0.1}\right) = 0 \\[1em]
& t_f < \infty
\end{array} \right.& 
\end{eqnarray*}
A few comments regarding the definition of Problem~$(Z_2)$ are in order:
\begin{enumerate}
\item Problem~$(Z_2)$ is a modification of Problem~$(Z_1)$ (which was a modification of Problem~$(Z_0)$). These problems are examples of the ancillary problems discussed in Section~\ref{sec:tyc2unscent}.
\item  The mean values of the final states are still required to be zero as in Problem~$(Z_1)$.  Thus, Problem~$(Z_2)$ levies additional conditions than Problem~$(Z_1)$.
\item Because the cost function is no longer the transit time, we set $t_f$ to be free as indicated by the statement $t_f < \infty$.
\item It is apparent that we cannot set $t_f \le t_f^b \approx 2.47$.  Consequently, we expect $t^2_{f} > 2.47$.  That is, we expect to pay a price for controlling the dispersion.
\end{enumerate}
The unscented instantiation of Problem~$(Z_2)$ is given by,
\begin{eqnarray*}
& \bu \in \U := \set{(u_1, u_2) \in \real{2}:\ u_1^2 + u^2_2 = 1 } & \nonumber\\
& \bX := \set{(x_i, y_i), i = 1, \ldots, N_\sigma} \in \real{2N_\sigma}   & \nonumber\\
&(Z_2^U) \left\{
\begin{array}{lll}
\text{Minimize }  & J_2^{U}[\bX(\cdot), \buf, t_f] :=\\
& \sum_{i=1}^{N_\sigma} w_i\, \left[\left(x_i(t_f)-\overline{x}(t_f)\right)^2 \right. + \\
& \left.\left.(y_i(t_f)-\overline{y}(t_f)\right)^2\right] \\
\text{Subject to} & \dot x_i(t) =  p_i y(t) + u_1(t) &\forall\ i =1, \ldots, N_\sigma\\
&\dot y_i(t) = q_i x(t) + u_2(t)  &\forall\ i =1, \ldots, N_\sigma\\
&(x_i(t_0), y_i(t_0)) = (2.25, 1)   &\forall\ i =1, \ldots, N_\sigma \\
&t_0 = 0 \\
&\overline{x}(t_f):= \sum_{i=1}^{N_\sigma} w_i\, x_i(t_f) = 0 \\
&\overline{y}(t_f):=\sum_{i=1}^{N_\sigma} w_i\, y_i(t_f) = 0 \\
\end{array} \right.& 
\end{eqnarray*}
As in \textbf{Steps~3} and \textbf{1}, Problem~$(Z_2^U)$ is a deterministic problem and can be solved by DIDO.

\subsection{Step 7: Solve and Analyze Problem~$(Z_2^U)$}
\begin{figure}[h!]
      \begin{center}
        \includegraphics[angle = 0,width = 0.48\columnwidth] {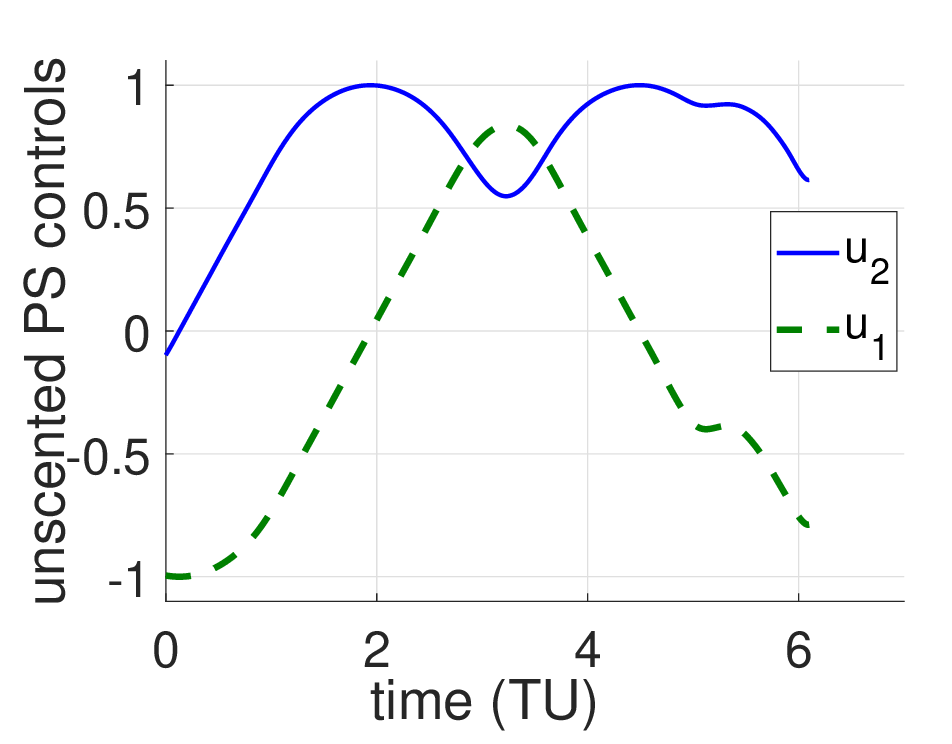}
        \includegraphics[angle = 0,width = 0.48\columnwidth] {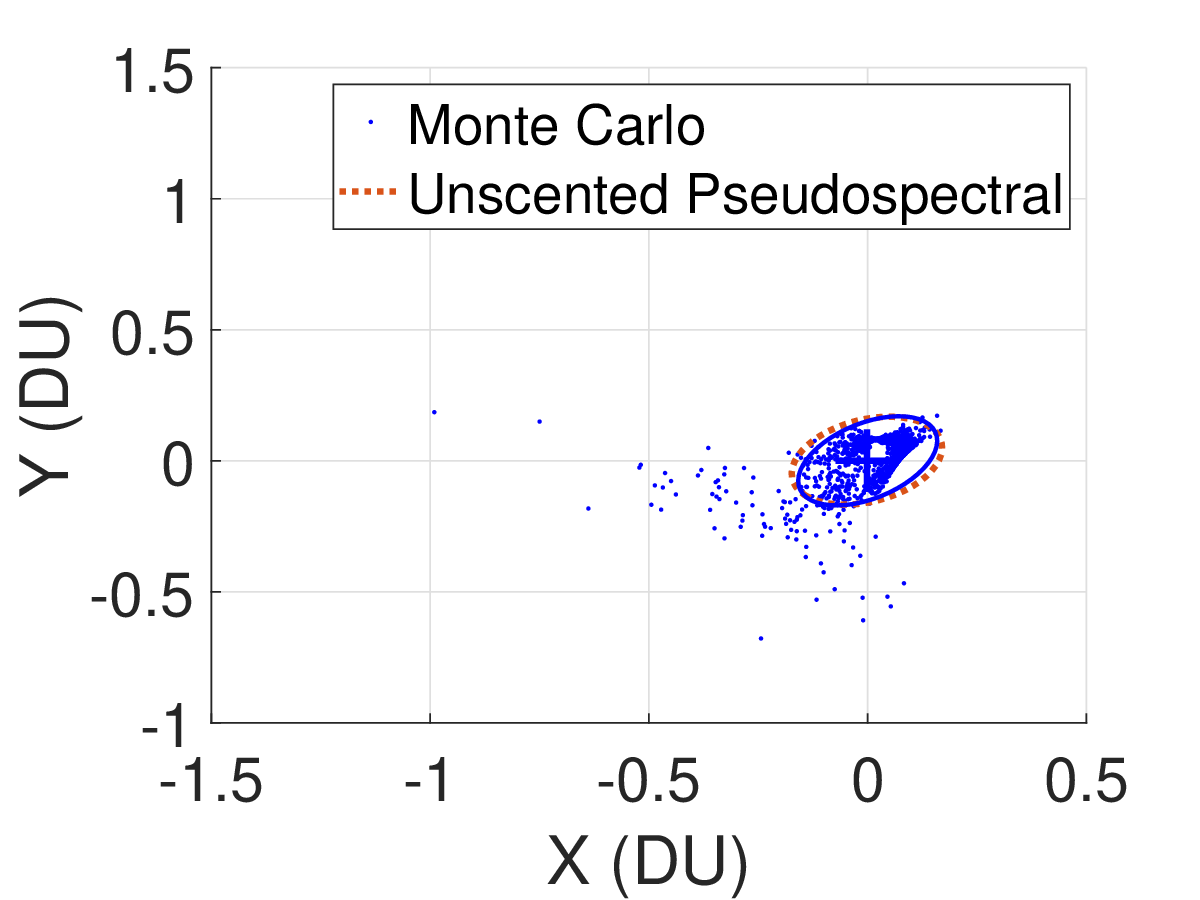}
      \caption{An unscented PS solution to Problem $Z_2$ and its Monte Carlo simulation.}
      \label{fig:Z2Sol}
      \end{center}
\end{figure}

The unscented PS control solution to Problem~$(Z_2^U)$ obtained via DIDO is shown in Fig.~\ref{fig:Z2Sol}.
It is clear that this unscented control solution is very different from the prior ones (compare Fig.~\ref{fig:Z1Sol}).
%
%
Post V\&V (see \textbf{Step~1}), this control was used to generate Monte Carlos simulations.  A result of this simulation is shown in Fig.~\ref{fig:Z2Sol}.  It is abundantly clear that the Monte Carlo simulations generate a very different and smaller spread of the endpoint values than those shown in Figs~\ref{fig:Z1Sol} or \ref{fig:MC4Z0}.

\begin{figure}[!h]
      \begin{center}
        \includegraphics[angle = 0,width = 0.75\columnwidth] {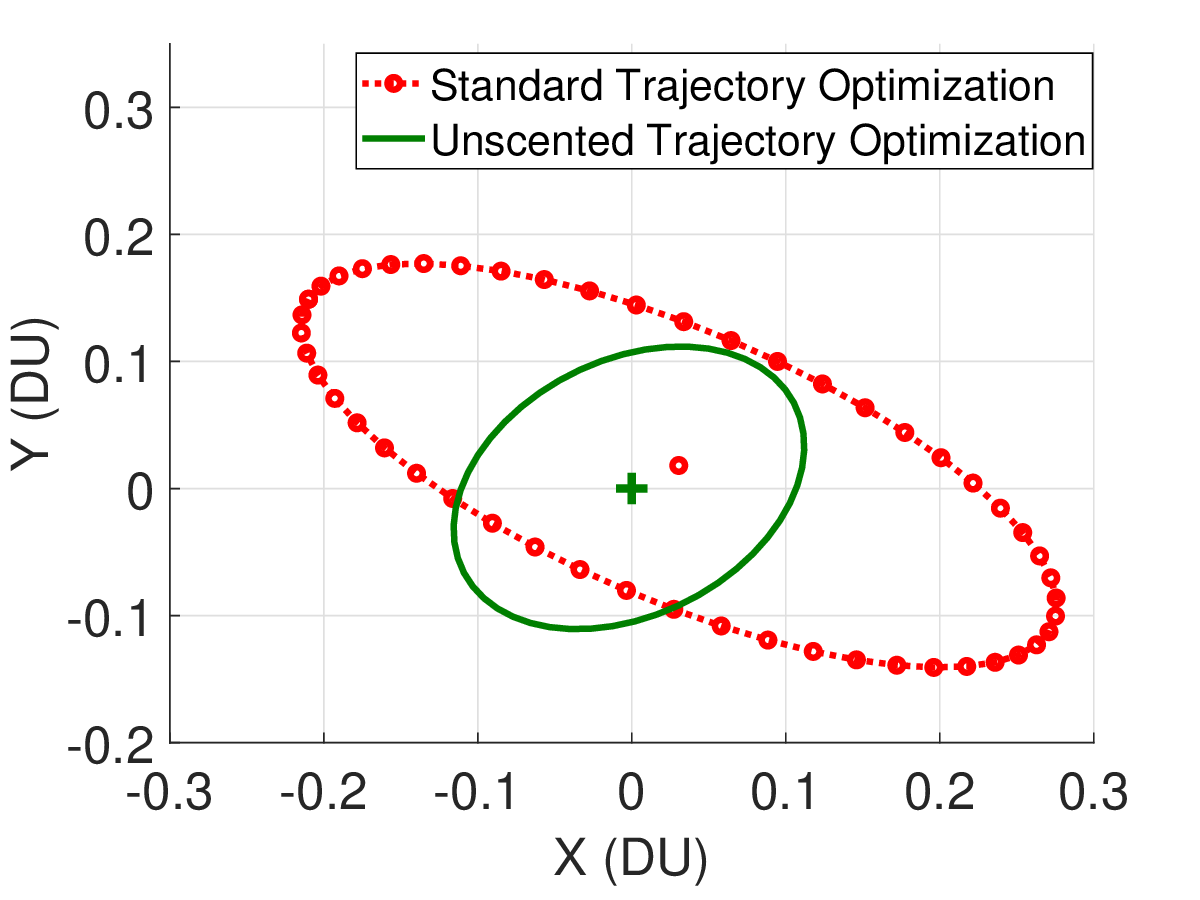}
      \caption{Covariance comparisons between standard (Problem~$(Z_0$)) and unscented (Problem~$(Z_2)$) solutions for uncertain wind conditions (given by \eqref{eq:p=PDFdata}).}
      \label{fig:MCCompare}
      \end{center}
\end{figure}

A direct comparison of covariance ellipses indicated in Figs.~\ref{fig:MC4Z0} and \ref{fig:Z2Sol} is shown in Fig.~\ref{fig:MCCompare}.
%
%
It is apparent that there is a significant improvement in the dispersions of the end point values. Note however that the transfer time is now significantly larger than $2.47$; it is, in fact approximately equal to $6.10$ time units (see the abscissa in Fig.~\ref{fig:Z2Sol}). In other words, unscented trajectory optimization has a price tag. Knowledge of this price tag during initial mission design may be critical to decision makers in exploring the trade space, particularly while making irreversible choices that have the potential to generate future requirements creep.

\subsection{Risk Estimates for Chance Constraints}
\begin{figure}[!h]
      \begin{center}
        \includegraphics[angle = 0,width = 0.8\columnwidth] {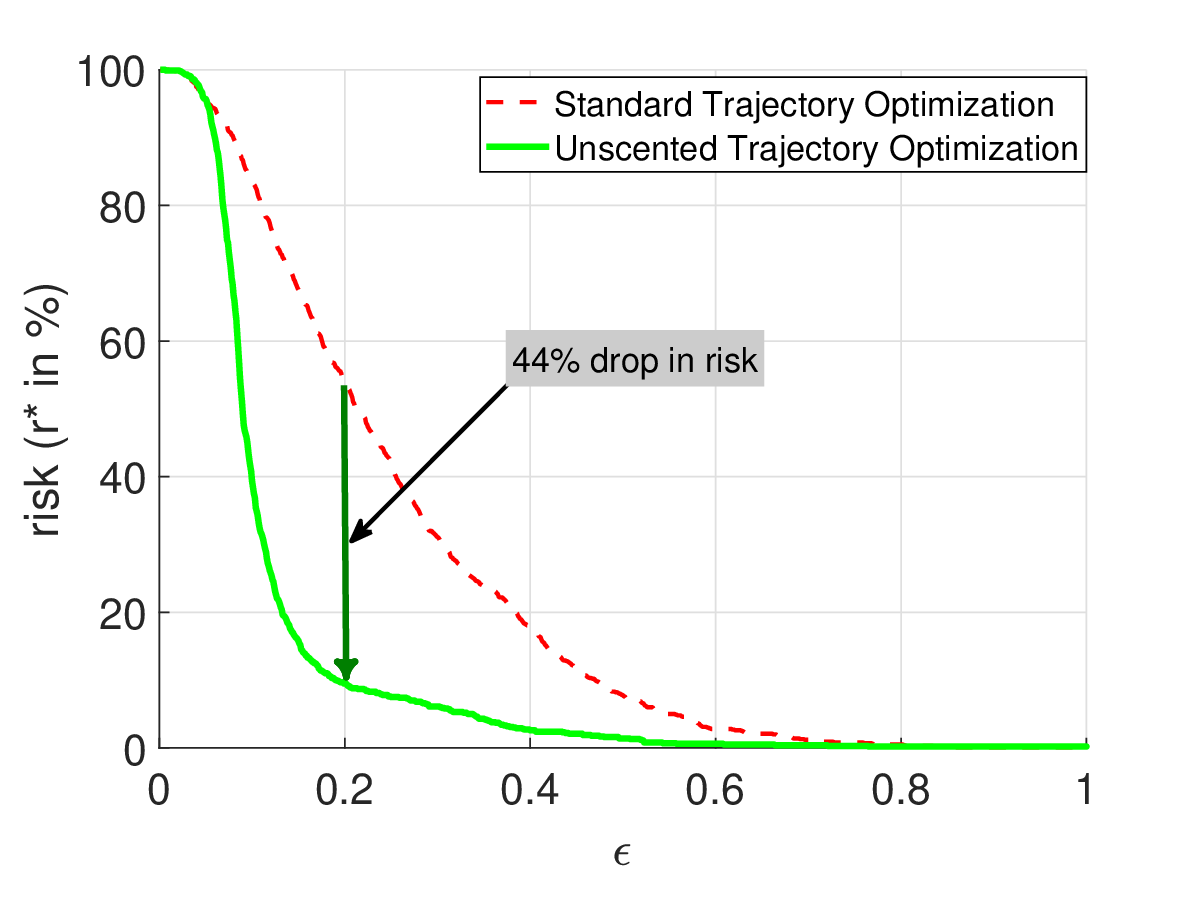}
      \caption{Risk estimates for a chance-constrained Zermelo problem.}
      \label{fig:risk}
      \end{center}
\end{figure}

Consider a chance-constrained Zermelo problem whose final-time condition is specified in the form of \eqref{eq:CC4T} where $\mathcal{T}$ is given by $\mathbb{B}(\bx^f, \varepsilon)$, a ball of radius $\varepsilon \ge 0$ centered at the point $\bx^f$.  In this case the chance constraint can be written as,
\begin{equation}\label{eq:chanceC4Z}
\text{Pr}\Big\{\bx(t_f, \bp) \in \mathbb{B}(\bx^f, \varepsilon) \Big\} \ge 1 - r
\end{equation}
where, $r \in (0, 1)$ is a specified level of risk. From \eqref{eq:Pr=040} it follows that if we choose $\varepsilon = 0$, then $r = 1$; i.e., the risk is $100\%$.  It thus follows that if $r$ were to be chosen to be any value less than one (with $\varepsilon = 0$), it is unlikely the chance-constrained Zermelo problem will even have a solution.  Hence the pair $(r, \varepsilon)$ must be chosen appropriately to ensure the mere existence of a solution to a chance-constrained problem. From the physics of the problem, it is apparent that for any given value of $\varepsilon$, there is a sufficiently large value of risk $r = r^*(\varepsilon)$ for which a solution exists. Estimates for $r^*(\varepsilon)$ can be easily obtained a posteriori from a Monte Carlo simulation.  These results are shown in Fig.~\ref{fig:risk} for Problems~$(Z_2)$ and $(Z_0)$.  The reduction in risk relative to the standard deterministic trajectory optimization is quite large.  For the sample case highlighted in Fig.~\ref{fig:risk}, for $\varepsilon = 0.2$, the $44\%$ drop in $r^*$ corresponds to about an $80\%$ reduction in risk relative to the standard optimal control. Note also that the drop in risk is quite steep for small increases in $\varepsilon$ near its smaller values.  For larger values of $\varepsilon$, the benefits are quite negligible.
%
%

\section{A Practical Example: Mission Recovery Under Subsystem Failure}
A number of real-world examples on the use of unscented optimal control techniques can be found in \cite{LS:acc-2,uoc-1:issfd,uoc-patent-2,Ng:UAV-2020,Karp:RS-experiment}. In this section we present the results obtained\cite{uoc-1:issfd} for recovering the mission of NASA's Hubble Space Telescope\cite{HST-point-1,HST-FGS,HST-moi,markley:HST-zero} (HST) using the techniques described in Section~\ref{sec:tyc2unscent}.

Gyros onboard the HST provide critical feedback signals for slewing the telescope while the fine guidance sensors (FGS) provide accuracy for precision pointing\cite{HST-point-1}. Normal operations of the HST requires four functioning gyros\cite{HST-point-1,markley:HST-zero}. Over the course of its mission, the HST gyros have failed and have been replaced by servicing missions.  With the cancelation of the servicing operations, NASA's challenge was to continue to operate the HST under failure modes\cite{markley:HST-zero}. A two-gyro mode of operation is possible\cite{TGS} but with reduced science productivity due to difficulties in slewing under constraints. A single gyro mode is also possible but at the price of a detrimental impact on the HST's science return\cite{cbsnews}.  The idea of unscented trajectory optimization for the HST was to provide zero-gyro attitude guidance signals to slew the spacecraft from a known point (made available by the FGS\cite{HST-FGS}) to near a target point for a handover to the FGS to provide terminal accuracy. Because the inertia dyadic of the HST is not known precisely\cite{HST-moi}, we consider the principal moments of inertia, $\bp:= (p_1, p_2, p_3) \in \real{3}$ to be the uncertain/tychastic variables and consider the nonlinear dynamical equations,
\begin{eqnarray}\label{eq:HST-dynamics}
\dot{q_1} & = & \frac{1}{2}\left[\omega_1q_4 - \omega_2q_3 + \omega_3q_2\right]\nonumber\\
\dot{q_2} & = & \frac{1}{2}\left[\omega_1q_3 + \omega_2q_4 - \omega_3q_1\right]\nonumber\\
\dot{q_3} & = & \frac{1}{2}\left[-\omega_1q_2 + \omega_2q_1 + \omega_3q_4\right]\nonumber\\
\dot{q_4} & = & \frac{1}{2}\left[-\omega_1q_1 - \omega_2q_2 - \omega_3q_3\right]\\
\dot{\omega_1} & = & \frac{u_1}{p_1} - \left(\frac{p_3 - p_2}{p_1}\right)\omega_2\omega_3\nonumber\\
\dot{\omega_2} & = & \frac{u_2}{p_2} - \left(\frac{p_1 - p_3}{p_2}\right)\omega_1\omega_3\nonumber\\
\dot{\omega_3} & = & \frac{u_3}{p_3} - \left(\frac{p_1 - p_2}{p_3}\right)\omega_1\omega_2\nonumber
\end{eqnarray}
where, $\bq := (q_1, q_2, q_3, q_4)$ and $\bomega:= (\omega_1, \omega_2, \omega_3) $ are the quaternions and body rates\cite{KaneLikinsLevinson}.  Note that not only are the dynamical equations nonlinear but so are also the dependencies of the dynamics function with respect to $\bp$.

\begin{figure}[!h]
      \begin{center}
      $\begin{array}{cc}
        \includegraphics[angle = 0,width = 0.48\columnwidth] {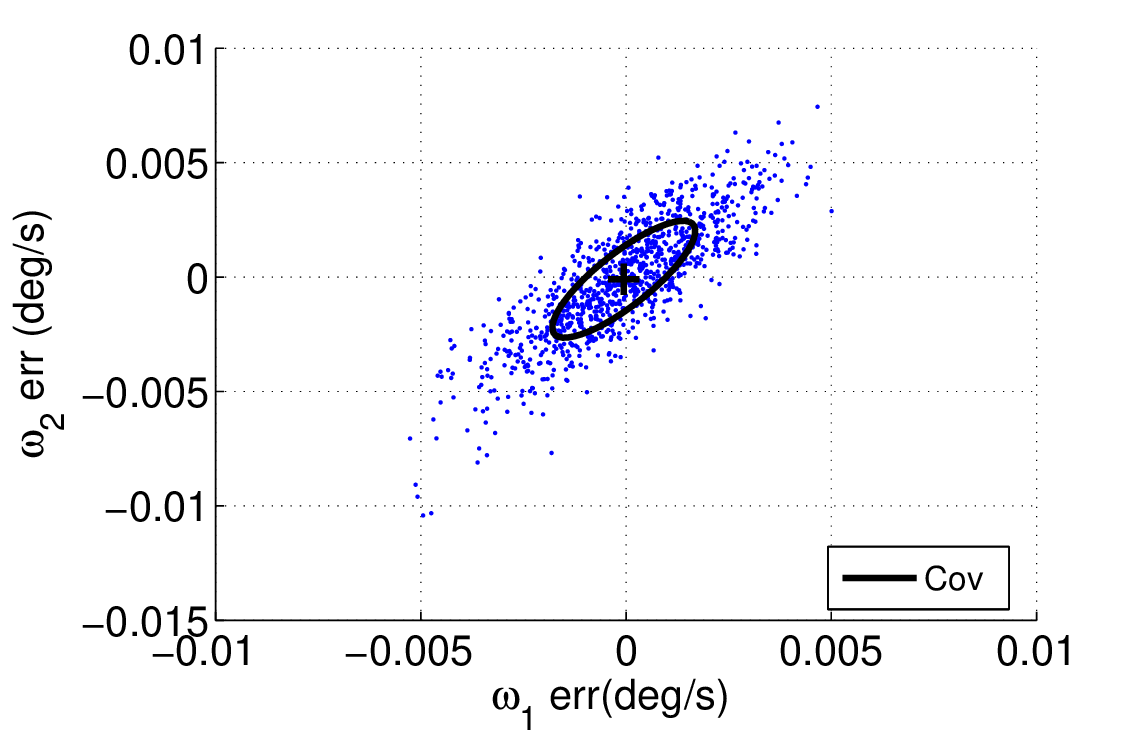}
        &
        \includegraphics[angle = 0,width = 0.48\columnwidth] {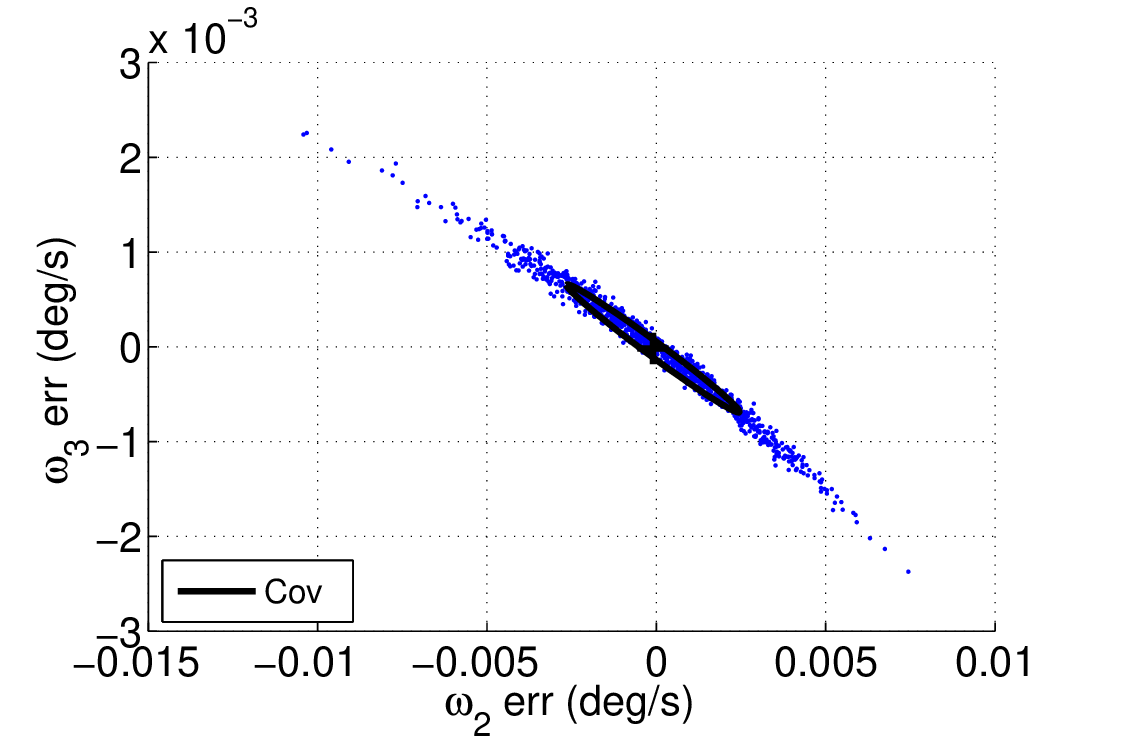}\\
         \includegraphics[angle = 0,width = 0.48\columnwidth] {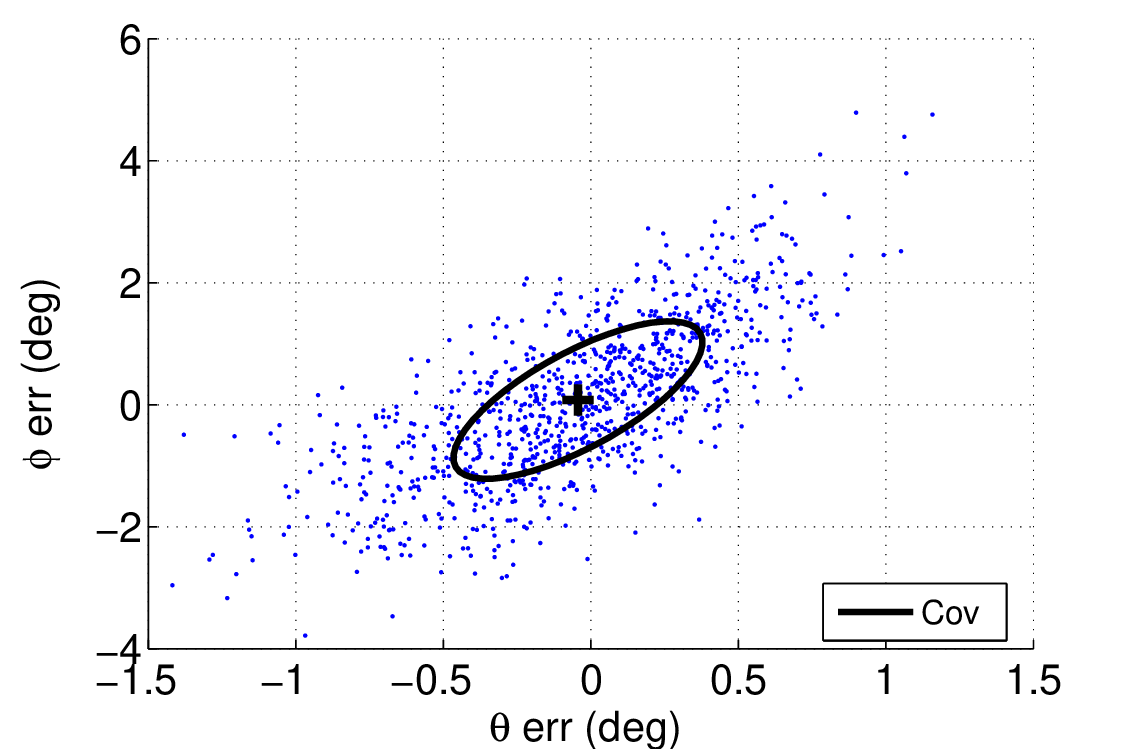}
        &
        \includegraphics[angle = 0,width = 0.48\columnwidth] {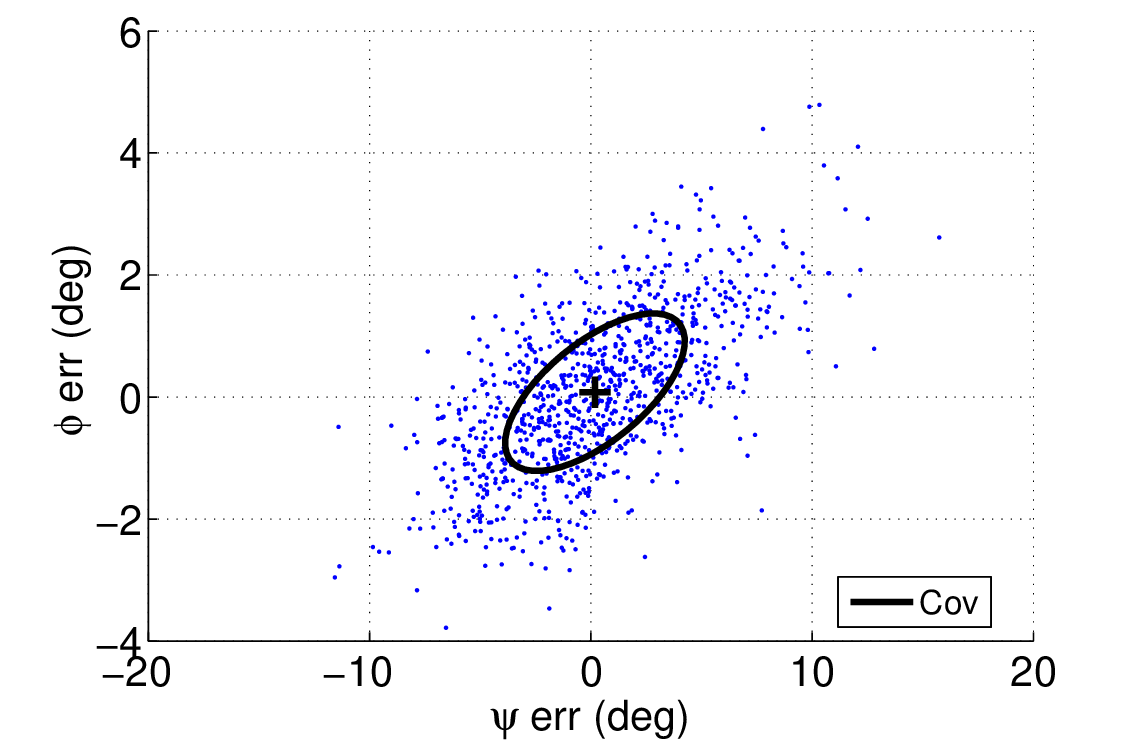}
      \end{array}$
      \caption{Monte Carlo simulations of target error distributions for the HST in a zero-gyro mode.}
      \label{fig:UQ-1}
      \end{center}
\end{figure}

\begin{figure}[!h]
      \begin{center}
      $\begin{array}{cc}
        \includegraphics[angle = 0,width = 0.48\columnwidth] {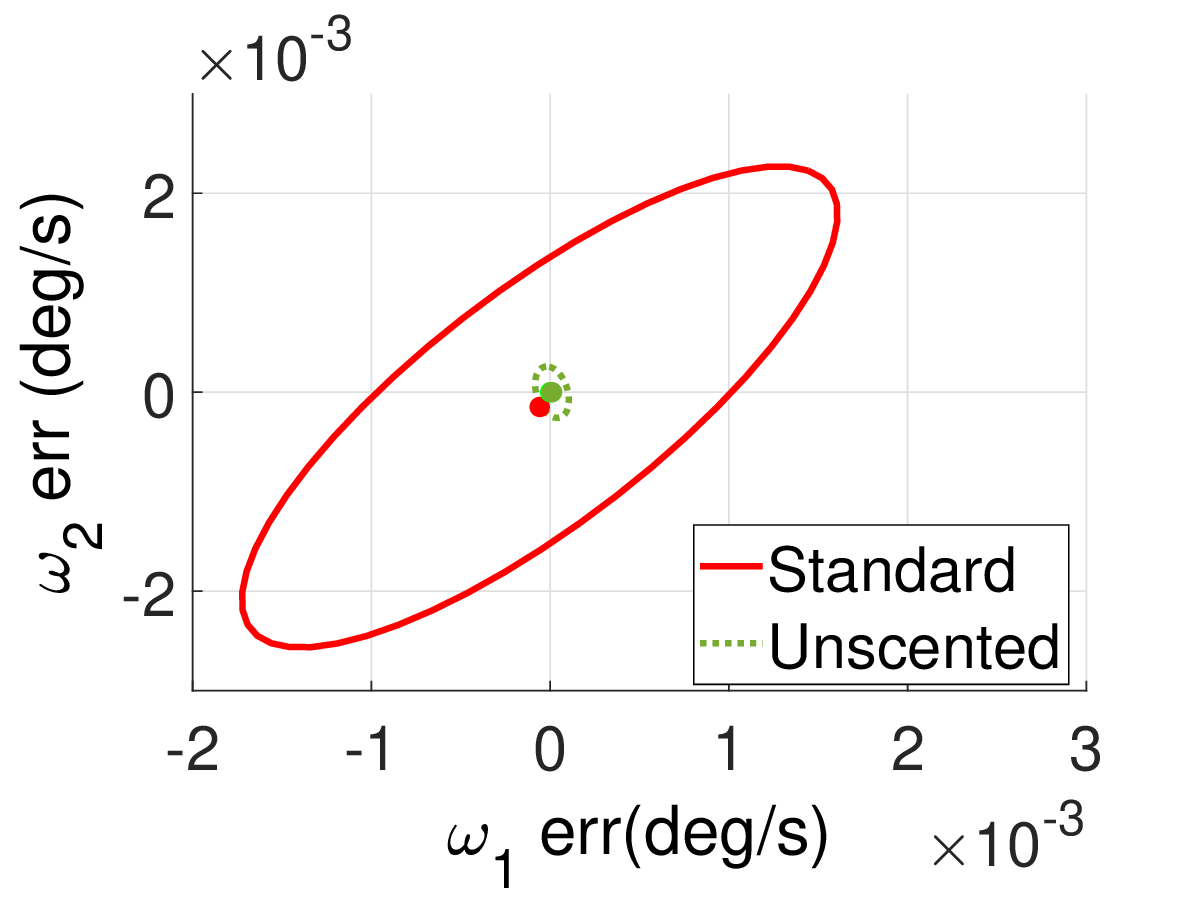}
        &
        \includegraphics[angle = 0,width = 0.48\columnwidth] {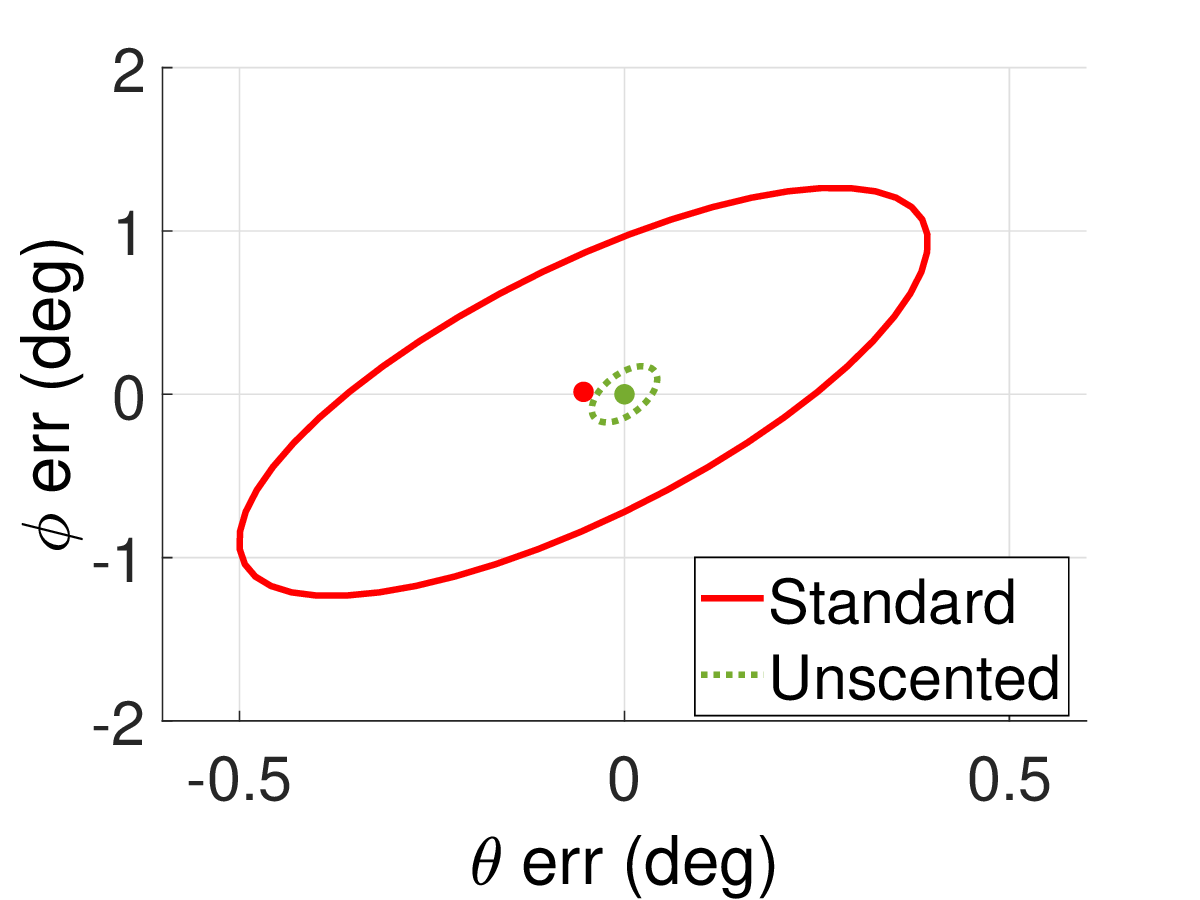}
      \end{array}$
      \caption{Performance of standard and unscented optimal control for the HST in a zero-gyro mode.}
      \label{fig:COVs-2}
      \end{center}
\end{figure}

For the baseline (deterministic) optimal control problem (see \textbf{Step~1}, in Section~\ref{sec:tyc2unscent}) we set the nominal value of $\bp =(36, 87, 94) \times 10^3 $ kg-m$^2$\cite{HST-moi} and consider the problem of slewing the HST in minimum time between two points given by,
%
\begin{eqnarray}\label{eq:boundary}
    \bx^0 &=& \left(0, 0, 0, 1, 0, 0, 0\right)\\
    \bx^f &=& \left(-0.27060, 0.27060, 0.65328, 0.65328, 0, 0, 0\right)\nonumber
\end{eqnarray}
%
These numbers corresponds to a yaw of $\psi(\bq) =90$ degrees, a pitch of $\theta(\bq) = 45$ degrees and a roll of $\phi(\bq) = 0$ degrees.  Fig.~\ref{fig:UQ-1} shows the errors in the targeted angular position and velocity generated by \textbf{Step~2}.  See \cite{uoc-1:issfd} for details of \textbf{Step~1}.
%
%
The plots shown in Fig.~\ref{fig:UQ-1} were obtained by considering a constrained Gaussian uncertainty of $3.3\%, 1 \sigma$ in the principal moments of inertia.  That is, the Gaussian variations in $\bp$ were constrained by the relevant physics constraints associated with an inertia dyad\cite{KaneLikinsLevinson}. Based on these results, it was concluded that the FGS would not be able to provide terminal guidance for the HST.  In order to recover the mission, the task of unscented guidance was to shrink the dispersions to a sufficiently small radius so that the FGS can take over the task of terminal guidance. Let,
\begin{equation}
\text{Cov}\left(\psi\big(\bq(t_f, \bp)\big), \theta\big(\bq(t_f, \bp)\big), \phi\big(\bq(t_f, \bp)\big), \omega_1(t_f, \bp), \omega_2(t_f, \bp), \omega_3(t_f, \bp)\right)
\end{equation}
denote the covariance matrix associated with the endpoint values of the random (tychastic) variables $\psi\big((\bq(t, \bp)\big)$, $\theta\big((\bq(t, \bp)\big), \phi\big((\bq(t, \bp)\big), \omega_1(t, \bp), \omega_2(t, \bp)$ and $\omega_3(t, \bp)$.  Note that $\psi, \theta$ and $\phi$ are not state variables.  To facilitate mission feasibility, the problem is to reduce the dispersions in the final values of $\psi, \theta, \phi, \omega_1, \omega_2$ and $\omega_3$ by constraining their variances according to,
\begin{multline}
\text{diag }\left[\text{Cov}\left(\psi\big(\bq(t_f, \bp)\big), \theta\big(\bq(t_f, \bp)\big), \phi\big(\bq(t_f, \bp)\big), \omega_1(t_f, \bp), \omega_2(t_f, \bp), \omega_3(t_f, \bp)\right)   \right] \le \\
(\sigma_{\psi}^2, \sigma_{\theta}^2,  \phi_{\phi}^2, \sigma_{\omega_1}^2, \sigma_{\omega_2}^2, \sigma_{\omega_3}^2)
\end{multline}
where $\sigma^2_{(\cdot)}$ are specified variances of the variables denoted by the subscript $(\cdot)$.
Carrying out the remainder of the steps delineated in Section~\ref{sec:tyc2unscent} up to the generation of Monte Carlo simulations, we get the results shown in  Table~\ref{table:uoc-risk}.
%
%
\begin{table}[h]
\center
\begin{tabular}
{l  c c  c r} \hline \hline

Target Parameter     & \multicolumn{2}{c}{Error Mean}       & \multicolumn{2}{c}{Variance}          \\
   & Unscented    & Standard    & Unscented   & Standard\\
\hline

 $ \omega_1 $(arcsec/s) & 0.036   & -0.197            & 0.341         &  6.26\\
 $ \omega_2 $(arcsec/s) & 0.000   &-0.356           & 0.932            & 9.18 \\
 $ \omega_3 $(arcsec/s) & 0.000  &-0.063            & 0.130           & 2.43 \\
 $ \psi $ (arcsec)     & 0.000   & 666            & 1,436          & 14,616\\
 $\theta$ (arcsec)     & 0.000     &-158            & 153            & 1,519\\
 $\phi$ (arcsec)       & 0.000   & 248            & 619            & 4,644\\

 \hline \hline
\end{tabular}
\caption{Comparison of Monte-Carlo-based error statistics between standard and unscented optimal control.}\label{table:uoc-risk}
\end{table}
%
%
Also shown in Table~\ref{table:uoc-risk} are the values of the corresponding statistics obtained using standard optimal control theory. From the values presented in Table~\ref{table:uoc-risk}, it is clear that the unscented solution not only targets the mean near perfectly but also that the variances in the key variables are reduced by an order of magnitude (see last two columns of Table~\ref{table:uoc-risk}). To better appreciate this latter point, the covariance ellipses between standard and unscented techniques are shown in Fig.~\ref{fig:COVs-2}.
%
%

\section{Conclusions}
The baseline risk associated with a solution to a static optimization problem with equality constraints is $100\%$. Consequently, probabilistic static optimization problems are framed for only inequalities. By definition, a trajectory optimization problem contains equalities in the form of dynamic constraints. Hence, a more careful consideration of uncertain variables appearing as equality constraints is necessary to formulate a probabilistic version of optimal control theory and its associated computational techniques.  Stochastic optimal control theory addresses this challenge in addition to the consideration of a Markov process. The theoretical challenge in using stochastic optimal control theory are the adjoint processes and the various Hamiltonians that involve the trace of an adjoint diffusion matrix and second-order adjoint equations.  Furthermore, a discretization of the stochastic differential equations involve intricate computational methods that account for the exotic rules of It\^{o} calculus.  Tychastic optimal control theory circumvents all these challenges. The tychastic perspective of unscented trajectory optimization ensures that the resulting computational problem formulation is consistent with the statistical mathematics associated with managing uncertainties across all functions and variables used to formulate a tychastic optimal control problem. The unscented transform and all its variants provide a fast and easy approach to generating a first-cut management of the proper statistics in tychastic trajectory optimization.  Together with sensitivity analysis, unscented trajectory optimization can be used for initial mission design to provide decision makers quick estimates of risk, reliability and confidence levels so that ``optimal missions'' do not suffer from cost overruns due to requirements creep.  In addition, the operational implementation of unscented controls can be used to enhance the safety of autonomous operations at distant celestial bodies particularly if/when there are failures in sensors or inner-loop feedback controls.

%
\section{Acknowledgment}
We thank the US Navy's Judge Advocate General's office for filing US Patents 9727034, 9849785, 10065312 and 10095198.  These patents describe further details of practical unscented trajectory optimization and their implementations.



\begin{thebibliography}{10}


\bibitem{kushner}
Kushner, H. J., ``Necessary Conditions for Continuous Parameter Stochastic Optimization Problems,'' \textit{SIAM Journal of Control}, Vol.~10, No.~3, 1972, pp.~550-565.

\bibitem{peng}
Peng, S., ``A General Stochastic Maximum Principle for Optimal Control Problems,'' \textit{SIAM Journal of Control and Optimization}, Vol.~28, No.~4, 1990, pp.~966-979.

\bibitem{yong-zhou-1999Book}
Yong, J. and Zhou, X. Y., \textit{Stochastic Controls: Hamiltonian Systems and HJB Equations}, Springer, New York, N.Y., 1999.

\bibitem{RMP-book}
Boltyanski, V. G. and Poznyak, A. S., \textit{The Robust Maximum Principle: Theory and Applications}, Birkh\"{a}user, New York, N.Y., 2012. 

\bibitem{SDE-sim-book}
Kloeden, P. E.,  Platen, E., \textit{Numerical Solution of Stochastic Differential Equations}, Springer, New York, N.Y., 1992.

\bibitem{SDE-num-2001}
Higham, D. J., ``An Algorithmic Introduction to Numerical Simulation of Stochastic Differential Equations,'' \textit{SIAM Review}, Vol.~43, No.3, 2001, pp.~525--546.


\bibitem{hnw-ode}
Hairer, E., N{\o}rsett, S. P. \& Wanner, G. \textit{Solving Ordinary Differential Equations I: Nonstiff Problems}, Springer-Verlag, 1993.

\bibitem{poznyak-2002}
Poznyak, A. S., ``Robust Stochastic Maximum Principle: Complete Proof and Discussions,'' \textit{Mathematical Problems in Engineering}, Vol.~8/4-5, 2002, pp.~389--411.


\bibitem{ross-book}
Ross, I. M., \textit{A Primer on Pontryagin's Principle in Optimal Control}, Second Edition, Collegiate Publishers, San Francisco, CA 2015.


\bibitem{shooting-2013}
Bonnans, J. F.,  ``The shooting approach to optimal control problems,''
\textit{IFAC Proceedings }, Volume 46, Issue 11, 2013, pp.~281--292.

\bibitem{trelat:survey}
Caillau, J.-B., Ferretti, R.,  Tr\'{e}lat, E.,  Zidani, H., ``An algorithmic guide for finite-dimensional optimal control problems,'' \textit{Handbook of Numerical Analysis}, Editor(s): E. Tr\'{e}lat, E. Zuazua, Volume 24,
2023, Elsevier, pp.~559--626.

\bibitem{conway:survey}
Conway, B. A., ``A Survey of Methods Available for the Numerical Optimization of Continuous Dynamic Systems,'' \textit{Journal of Optimization Theory and Applications}, Vol.~152, 2012, pp.~271--306.


\bibitem{uoc-1:issfd}
Ross,  I. M., Proulx R. J. and M. Karpenko, ``Unscented Optimal Control For Space Flight,'' \textit{24th International Symposium on Space Flight Dynamics (ISSFD)}, Laurel, MD, May 5-9, 2014.

\bibitem{uoc-2:SD}
Ross, I. M., Proulx, R. J., and Karpenko, M., ``Unscented Optimal Control for Orbital and Proximity Operations in an Uncertain Environment: A New Zermelo Problem,'' \textit{AIAA Space and Astro. Forum and Expo.: AIAA/AAS Astrodynamics Specialist Conf.}, San Diego, CA 4-7 August 2014.

\bibitem{ug:acc}
Ross, I. M., Proulx, R. J., and Karpenko, M. ``Unscented Guidance,'' \textit{Proceedings of the ACC}, Chicago, IL, July 1--3, 2015.

\bibitem{RS:jgcd}
Ross, I. M., Proulx, R. J.  Karpenko, M. and Q. Gong, ``Riemann-Stieltjes Optimal Control Problems for Uncertain Dynamcial Systems,'' \textit{Journal of Guidance, Control and Dynamics}, Vol.~38, No.~7, 2015, pp.~1251--1263.

\bibitem{LS:acc-1}
Ross, I. M., Karpenko, M. and Proulx, R. J., ``A Lebesgue-Stieltjes Framework For Optimal Control and Allocation,'' \textit{Proceedings of the ACC}, Chicago, IL, July 1--3, 2015.

\bibitem{uo-1:2015}
Ross, I. M., Proulx, R. J. and Karpenko, M., ``Unscented Optimization,'' \textit{AAS/AIAA Astrodynamics Specialist Conference}, Vail, CO, August 9--13, 2015. AAS Paper 15-607.

\bibitem{LS:acc-2}
Ross, I. M., Karpenko, M. and Proulx, R. J., ``Path Constraints in Tychastic and Unscented Optimal Control: Theory, Application and Experimental Results,'' \textit{Proceedings of the ACC}, Boston, IL, July 6--8, 2016.

\bibitem{uoc-patent-1}
Ross, I. M., Proulx, R. J., Karpenko, M ``Unscented Control for Uncertain Dynamic Systems,'' US Patent 9727034B1, April 29, 2015.


\bibitem{uo-1:MRYpts-2016}
Ross, I. M., Proulx, R. J. and Karpenko, M., ``Monte Rey Methods for Unscented Optimization,'' \textit{AIAA Guidance, Navigation, and Control Conference}, AIAA SciTech Forum, San Diego, CA, January 4--8, 2016. doi:10.2514/6.2016-0871.


\bibitem{uoc-patent-2}
Ross, I. M., Karpenko, M. and Proulx, R. J., ``Method and Apparatus for State Space Trajectory Control of Uncertain Dynamic Systems,'' US Patent 9849785B1, June 21, 2016.

\bibitem{uo-patent-2016}
Ross, I. M., Karpenko, M. and Proulx, R. J., ``Unscented Optimization and Control Allocation,'' US Patent 10065312B1, July 13, 2016.

\bibitem{uo-patent-2017}
Ross, I. M., Karpenko, M. and Proulx, R. J., ``Closed-Loop Control System Using Unscented Optimization,'' US Patent 10095198B1, July 06, 2017.



\bibitem{aubin=tyche}
Aubin, J.-P., Chen, L., Dordan, O., Faleh, A., Lezan, G. and Planchet,
F., ``Stochastic and Tychastic Approaches to Guaranteed ALM Problem,''
\textit{Bulletin Fran\c{c}ais d'Actuariat}, Vol. 12, No. 23, 2012, pp. 59-95.

\bibitem{clsw}
Clarke, F. H., Ledyaev, Y. S.,  Stern, R. J., and Wolenski, P. R.,
{\it Nonsmooth Analysis and Control Theory,} Springer-Verlag, New
York, NY, 1998.

\bibitem{koopman-2} Koopman, B. O., ``The Theory of Search, II: Target Detection,'' \textit{Operations Research}, Vol.~4, No.~5, 1956, pp.~503--531.



\bibitem{stone}
Stone, L. D. and Richardson, H. R., ``Search for Targets with Conditionally Deterministic Motion,'' \textit{SIAM Journal of Applied Mathematics}, Vol.~27, No.~2, September 1974, pp.~239--255.

\bibitem{phelps-2:cdc}
Phelps, C., Royset, J. O., and Gong, Q., ``Sample Average Approximations
in Optimal Control of Uncertain Systems,'' \textit{Proceedings of the 52nd IEEE
CDC}, Florence, Italy, 2013.

\bibitem{Li:phys}
Li, J.-S. and Khaneja, N., ``Control of Inhomogeneous Quantum Ensembles,'' \textit{Physical Review A}, Vol.~73, No.~030302(R), 2006, pp.~030302-1--030303-4.

\bibitem{Li:pnas}
Li, Jr-S., Ruths, J., Yu, T-Y,  Arthanari H., and Wagner, G., ``Optimal Pulse Design in Quantum Control: A Unified Computational Method,'' \textit{Proceedings of the NAS}, Vol.~108, No.~5, Feb 2011, pp.~1879-1884.

\bibitem{RussiaReport}
Kurzhanski, A. B. (Ed), \textit{Advances in Nonlinear Dynamics and Control: A Report from Russia}, Birkh\"{a}user, Boston, 1993.


\bibitem{vinter-minimax}
Vinter, R. B., ``Minimax Optimal Control,'' \textit{SIAM Journal of Control
and Optimization}, Vol. 44, No. 3, 2005, pp. 939-968.


\bibitem{aubin-alm-2012}
Aubin, J.-P., Chen, L., Dordan, O., Faleh, A., Lezan, G. and Planchet, F., ``Stochastic and Tychastic Approaches to Guaranteed ALM Problem,'' \textit{Bulletin Fran{\c c}ais d'Actuariat}, Vol.~12, No.~23, 2012, pp.~59-95.


\bibitem{julier:simplex}
Julier, S., J., ``The Spherical Simplex Unscented Transformation,'' \textit{Proceedings of the American Control Conference}, Denver, CO, June 4-6, 2003, Vol.~3, pp.~2430--2434.

\bibitem{julier-acc-95}
Julier, S. J., Uhlmann, J. K. and H. F. Durrant-Whyte, ``A New Approach for Filtering Nonlinear Systems,'' \textit{Proceedings of the  American  Control Conference}, Seattle, WA, 1995, pp.~1628--1632.


\bibitem{longuski}
Longuski, J. M., Guzm\'{a}n, J, J., and Prussing, J. E., \textit{Optimal Control with Aerospace Applications}, Springer, New York, N.Y., 2014.

\bibitem{bryson:ho}
Bryson, A.E., and Ho, Y.C., {\it Applied Optimal Control,}
Hemisphere, New York, 1975.


\bibitem{vinter}
 Vinter, R. B. \textit{Optimal Control}, Birkh\"{a}user, Boston, MA, 2000.



\bibitem{charnes-cooper}
Charnes A. and  Cooper, W. W., ``Chance-Constrained Programming,'' \textit{Management Science}, Vol.~6, No.~1, 1959, pp.~73--79.

\bibitem{engels}
Engels, H., \textit{Numerical Quadrature and Cubature}, Academic Press, New York, 1980.

\bibitem{cools_2001}
Cools, R., Mysovskikh, I. P., and Schmid, H. J., ``Cubature Formulae and Orthogonal Polynomials,'' \textit{Journal of Computational and Applied Mathematics}, Vol.~127, pp.~121-152, 2001.


\bibitem{UT-various-2012}
Adurthi, N., Singla, P. and Singh, T.,  ``The Conjugate Unscented Transform - An Approach to Evaluate Multi-dimensional Expectation Integrals,'' \textit{Proceedings of the American Control Conference}, Montreal, QC, Canada, 2012, pp.~5556--5561.


\bibitem{UT-various-2015}
Menegaz, H. M. T.,  Ishihara, J. Y.,  Borges, G. A. and Vargas, A. N., ``A Systematization of the Unscented Kalman Filter Theory,'' \textit{IEEE Transactions on Automatic Control}, vol. 60, no. 10, pp. 2583--2598, 2015.

\bibitem{UT-various-2021}
Stojanovski, Z. and Savransky, D., ``Higher-Order Unscented Estimator,'' \textit{Journal of Guidance, Control, and Dynamics}, 44/12, 2021, 2186--2198.


\bibitem{DIDO:arXiv}
Ross, I. M., ``Enhancements to the DIDO Optimal Control Toolbox,'' arXiv preprint, arXiv:2004.13112, 2020, https://arxiv.org/abs/2004.13112.

\bibitem{RossKarp_IFAC_2012}
Ross, I.~M. and Karpenko, M.
``A Review of Pseudospectral Optimal Control: From Theory to
  Flight,''  \textit{Annual Reviews in Control}, Vol.~36, No.~2, pp. 182--197, 2012.

\bibitem{spec-alg}
Gong, Q., Fahroo F. and Ross, I. M., ``Spectral Algorithm for Pseudospectral Methods in Optimal Control,'' \textit{Journal of Guidance, Control, and Dynamics}, vol. 31 no. 3, pp. 460-471, 2008.


\bibitem{auto-knots}
Gong, Q. and  Ross, I. M., ``Autonomous Pseudospectral Knotting Methods for Space Mission Optimization,'' \textit{Advances in the Astronatuical Sciences}, Vol.~124, 2006, AAS 06-151, pp.~779--794.


 \bibitem{arb-grid}
Gong, Q., Ross, I. M. and Fahroo, F., ``Spectral and Pseudospectral Optimal Control Over Arbitrary Grids,'' \textit{Journal of Optimization Theory and Applications}, vol.~169, no.~3, pp.~759-783, 2016.


\bibitem{fastmesh}
Koeppen, K., Ross, I. M.,  Wilcox, L. C. and Proulx, R. J., ``Fast Mesh Refinement in Pseudospectral Optimal Control,'' \textit{Journal of Guidance, Control, and Dynamics}, vol. 42 no. 4, pp. 711-722, 2019.

\bibitem{BirkhoffTN}
Ross, I. M. and Proulx, R. J., ``Further Results on Fast Birkhoff Pseudospectral Optimal Control Programming,'' \textit{Journal of Guidance, Control, and Dynamics}, vol. 42 no. 9, pp. 2086-2092, 2019.

\bibitem{IEEE:spectrum}
Bedrossian, N., Karpenko, M., and Bhatt, S.,
``Overclock My Satellite: Sophisticated Algorithms Boost Satellite Performance on the Cheap,''
\textit{IEEE Spectrum}, November 2012.


\bibitem{Bedrossian_2009}
Bedrossian, N., Bhatt, S., Kang, W., and Ross, I.~M., ``Zero-Propellant
  Maneuver Guidance,'' \textit{IEEE Control Systems Magazine}, October 2009,
  pp.~53--73.


\bibitem{SIAMnews}
Kang, W. and Bedrossian, N., ``Pseudospectral Optimal Control Theory
Makes Debut Flight -- Saves NASA \$1M in Under 3 hrs,'' \textit{SIAM
News}, Vol.~40, No.~7, September 2007, Page 1.

\bibitem{LRO-CSM}
Karpenko, M. et al, ``Fast Attitude Maneuvers for NASA's Lunar Reconnaissance Orbiter: Practical Flight Application of Attitude Guidance Using Birkhoff Pseudospectral Theory and Hamiltonian Programming,''
\textit{IEEE Control Systems Magazine}, to appear.


\bibitem{serres}
Serres, U., ``On Zermelo-Like Problems: Gauss-Bonnet Inequality and E. Hopf Theorem,'' \textit{Journal of Dynamical and Control Systems}, January 2009, Vol. 15, No. 1, pp.~99-131.


\bibitem{kamel}
Kamel, A., and Tibbitts, R., ``Some Useful Results on Initial Node Locations for Near-Equatorial Circular Satellite Orbits,'' \textit{Celestial Mechanics}, Vol. 8, 1973, pp. 45-73.

\bibitem{kech}
Kechichian, J. A., ``Optimal Steering for North-South Stationkeeping of Geostationary Spacecraft,'' \textit{Journal Of Guidance, Control, And Dynamics}
Vol.~20, No.~3, May-June 1997, pp.~435-444.


\bibitem{stevens-Zerm}
Stevens, R. E. and Baker, W. P., ``Optimal Control of a Librating Electrodynamic Tether Performing a Multirevolution Orbit Change,'' \textit{AIAA Guidance, Navigation, and Control Conference} 10-13 August 2009, Chicago, Illinois. AIAA 2009-5803.


\bibitem{Karp:RS-experiment}
Karpenko, M. and Proulx, R. J., ``Experimental Implementation of Riemann-Stieltjes Optimal Control for Agile Imaging Satellites,'' \textit{Journal of Guidance, Control, and Dynamics}, 39/1, 2016, pp.~144-150.

\bibitem{Ng:UAV-2020}
Ng, H. K., ``Strategic Planning with Unscented Optimal Guidance for Urban Air Mobility,'' \textit{AIAA Aviation Forum}, 2020, https://doi.org/10.2514/6.2020-2904



\bibitem{HST-point-1}
Dougherty, H., Tompetrini, K. Levinthal, J. and Nurre, G., ``Space Telescope Pointing Control System,'' \textit{Journal of Guidance, Control and Dynamics}, Vol.~5, No.~4, 1982, pp.~403-409.

\bibitem{HST-FGS}
Beals, G. A., Crum, R. C., Dougherty, H. J., Hegel, D. K., Kelley, J. L., and Rodden, J. J., ``Hubble Space Telescope Precision Pointing Control System,'' \textit{Journal of Guidance, Control and Dynamics}, Vol.~11, No.~2, 1988, pp.~119-123.

\bibitem{HST-moi}
Thienel, J. K., and Sanner, R. M., ``Hubble Space Telescope Angular Velocity Estimation During the Robotic Servicing Mission,'' \textit{Journal of Guidance, Control, and Dynamics}, Vol.~30, No.~1, 2007, pp.~29-34.


\bibitem{markley:HST-zero}
Markley, F. L. and Nelson, ``Zero-Gyro Safemode Controller for the Hubble Space Telescope,''
\textit{Journal of Guidance, Control and Dynamics}, Vol.~17, No.~4, 1994, pp.~815-822.

\bibitem{TGS}
Prior, M. and Dunham, L., ``System Design and Performance of Two-Gyro Scince Mode for the Hubble Space Telescope,'' \textit{Acta Astronautica}, Vol.~61, 2007, pp.~1010-1018.

\bibitem{cbsnews}
Harwood, W., ``Healthy Hubble Telescope Raises Hopes of Longer Life,'' \textit{Spaceflightnow.com}, Posted June 1, 2013.

\bibitem{KaneLikinsLevinson}
Kane, T. R., Likins, P. W. and Levinson, D. A., \textit{Spacecraft Dynamics}, McGraw-Hill, 1983.











\end{thebibliography}

\end{document}